\documentclass[12pt, reqno]{amsart}
\usepackage{hyperref}
\usepackage{amssymb, amsthm, amsmath, amsfonts}
\usepackage{mathrsfs} 
\usepackage{array, epsfig}
\usepackage{bbm}
\usepackage[numbers,square]{natbib}

  \evensidemargin
\oddsidemargin

\newcommand{\bA}{\mathbb{A}}

\newcommand{\bB}{{\bar {\mathbb{B}}}}

\newcommand{\bP}{\mathbb{P}}

\newcommand{\bT}{\mathbb{T}}

\newcommand{\cN}{\mathcal{N}}

\newcommand{\Normal}{\mathcal{N}}

\newcommand{\E}{\mathbb E}

\newcommand{\R}{\mathbb{R}}
\newcommand{\N}{\mathbb{N}}

\newcommand{\C}{\mathbb{C}}

\renewcommand{\P}{\mathbb{P}}

\renewcommand{\Re}{\operatorname{Re}}
\renewcommand{\Im}{\operatorname{Im}}

\newcommand{\Var}{\mathop{\mathrm{Var}}\nolimits}

\newcommand{\Cov}{\mathop{\mathrm{Cov}}\nolimits}

\newcommand{\eps}{\varepsilon}

\newcommand{\todistr}{\overset{d}{\underset{n\to\infty}\longrightarrow}}

\newcommand{\toweak}{\overset{w}{\underset{n\to\infty}\longrightarrow}}
\newcommand{\toprobab}{\overset{P}{\underset{n\to\infty}\longrightarrow}}
\newcommand{\toprobabnon}{\overset{P}{\longrightarrow}}

\newcommand{\ton}{\overset{}{\underset{n\to\infty}\longrightarrow}}

\newcommand{\ind}{\mathbbm{1}}

\newcommand{\dd}{{\rm d}}
\newcommand{\eee}{{\rm e}}

\newcommand{\limto}[2] {\lim_{#1 \to #2}}
\newcommand{\indi}[1]{\mathbbm{1}_{ \left\{ #1 \right\} }}

\newcommand{\bS}{\mathbf{S}}
\newcommand{\bZ}{\mathbf{Z}}

\theoremstyle{plain}
\newtheorem{theorem}{Theorem}[section]
\newtheorem{lemma}[theorem]{Lemma}
\newtheorem{corollary}[theorem]{Corollary}
\newtheorem{proposition}[theorem]{Proposition}

\theoremstyle{definition}

\newtheorem{example}[theorem]{Example}

\newtheorem{conjecture}[theorem]{Conjecture}

\theoremstyle{remark}
\newtheorem{remark}[theorem]{Remark}

\begin{document}

\author{Zakhar Kabluchko}
\address{Zakhar Kabluchko: Institut f\"ur Mathematische Stochastik,
Westf\"alische Wilhelms-Universit\"at M\"unster,
Orl\'eans--Ring 10,
48149 M\"unster, Germany}
\email{zakhar.kabluchko@uni-muenster.de}

\author{Hauke Seidel}
\address{Hauke Seidel: Institut f\"ur Mathematische Stochastik,
Westf\"alische Wilhelms-Universit\"at M\"unster,
Orl\'eans--Ring 10,
48149 M\"unster, Germany}
\email{hauke.seidel@uni-muenster.de}

\email{}

\title[Critical points of random polynomials]{Distances between zeroes and critical points for random polynomials with i.i.d.\ zeroes}


\thanks{}


\begin{abstract}
Consider a random polynomial $Q_n$ of degree $n+1$ whose zeroes are i.i.d.\ random variables $\xi_0,\xi_1,\ldots,\xi_n$ in the complex plane. We study the pairing between the zeroes of $Q_n$ and its critical points, i.e.\  the zeroes of its derivative $Q_n'$. In the asymptotic regime when $n\to\infty$, with high probability there is a critical point of $Q_n$ which is very close to $\xi_0$. We localize the position of this critical point by proving that the difference between $\xi_0$ and the critical point has approximately complex Gaussian distribution with mean $1/(nf(\xi_0))$ and variance of order $\log n \cdot  n^{-3}$.
Here, $f(z)= \E[\frac 1 {z-\xi_k}]$ is the Cauchy--Stieltjes transform of the $\xi_k$'s. We also state some conjectures on critical points of polynomials with dependent zeroes, for example the Weyl polynomials and characteristic polynomials of random matrices.
\end{abstract}

\maketitle

\section{Introduction}
Critical points of a polynomial $Q$ are defined as complex zeroes of its derivative $Q'$.
The Gauss--Lucas theorem states that the critical points of any polynomial are contained in the convex hull of its zeroes. Numerous results on the location of the zeroes and the critical points of deterministic polynomials can be found in the book~\cite{rahman_schmeisser_book}.

In this paper, we shall be interested in  random polynomials. Let $\xi_0,\xi_1,\ldots$ be a sequence of independent and identically distributed random variables taking complex values.
Consider the random polynomial
$$
Q_n(z) := \prod_{k=0}^n (z - \xi_k).
$$
The study of critical points of such polynomials was initiated by Pemantle and Rivin~\cite{pemantle_rivin}. Confirming their conjecture, one of the authors proved in~\cite{kabluchko15} that that the empirical probability measure
\begin{align}
\mu_n = \frac{1}{n} \sum_{z\in\C\colon Q_n'(z)=0} \delta_z \label{Kabluchko15}
\end{align}
counting (with multiplicities) the critical points of $Q_n$ converges in probability (and weakly) to the probability distribution of $\xi_0$. That is, for large $n$, the critical points have approximately the same distribution as the zeroes.
The reader should keep in mind that this does not necessarily hold in the deterministic setting. For example, the zeros of the polynomial $Q(z)=z^n-1$ are the $n$-th roots of unity, which are all on the unit circle, while the critical points of $Q$ are all equal $0$. Further results in this direction were obtained in~\cite{subramanian,subramanian_phd,orourke_matrices,orourke,hu_chang,byun,Tulasi16,tulasi_phd}.

\begin{figure}[!tbp]
\centering
\includegraphics[width=0.75\textwidth]{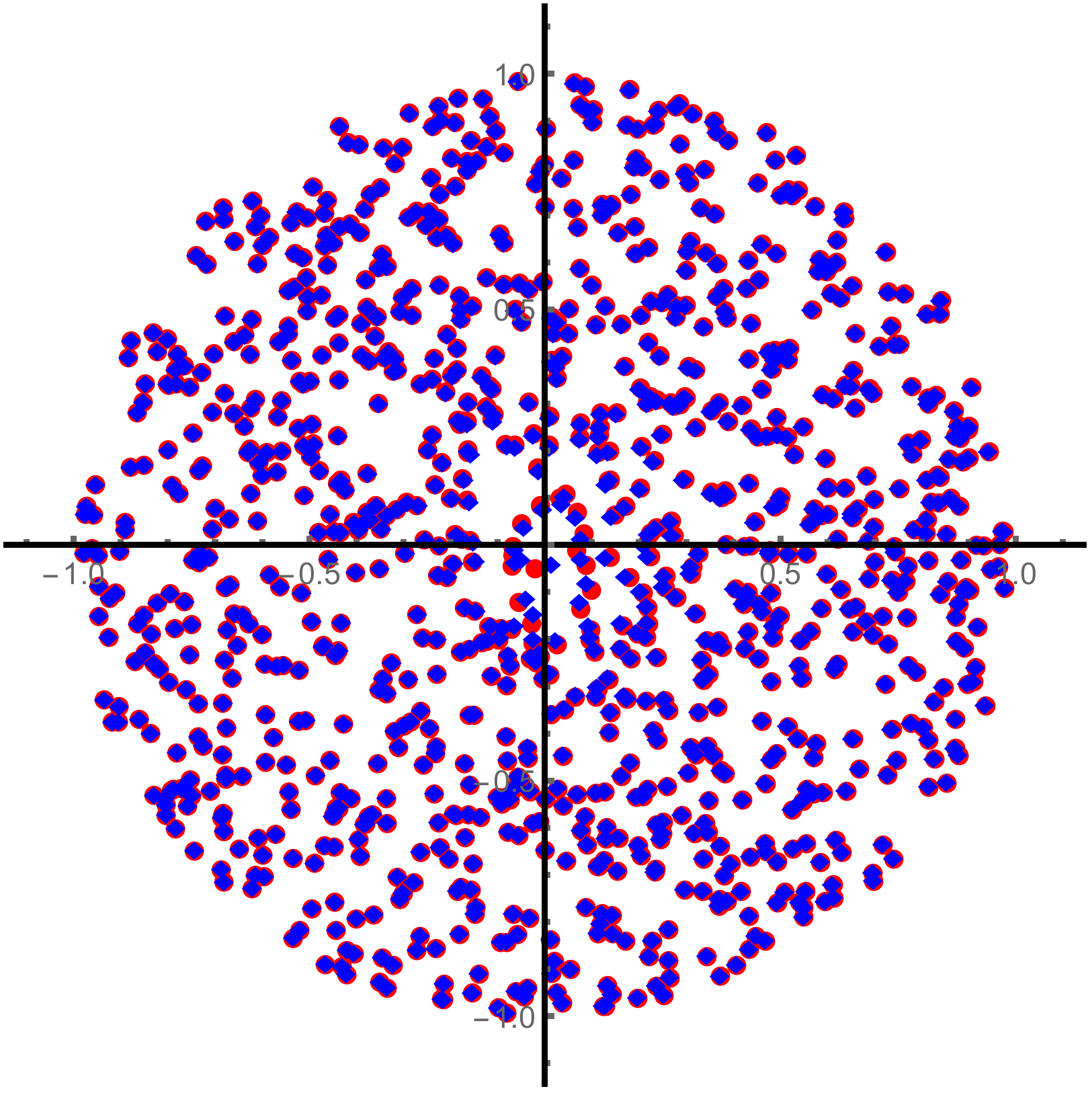}
\caption{Zeroes and critical points of a random polynomial of degree $n=1000$ whose i.i.d.\ zeroes have uniform distribution on the unit disk. Red disks: zeroes. Blue diamonds: critical points.}
	\label{fig:zeroes_critical}
\end{figure}

It was observed by Hanin~\cite{hanin_riemann,hanin_gauss,hanin_poly} that the zeroes and critical points of various random large-degree polynomials tend to appear in pairs; see Figure~\ref{fig:zeroes_critical}. More precisely,  the distance between the zero and the closest critical point is usually much smaller than the typical distance between close zeroes. Among other results, Hanin~\cite{hanin_poly} localised the position of the critical point associated to some fixed zero up to an error term of order $o(1/n)$.
The purpose of the present article is to prove a central limit theorem describing the random fluctuations of the critical point near its expected position.
We deal with polynomials having i.i.d.\ zeroes, as defined above, but numerical simulations suggest that some of the results should hold for other ensembles of random polynomials. This will be discussed in Section~\ref{sec:conjectures}.

\subsection{Notation}
Let $\bB_{r}(u)= \{z\in\C\colon |z-u| \leq  r\}$ be the closed disk of radius $r$ centered at $u\in\C$. Let $\cN_\C(0,\sigma^2)$ denote a complex normal distribution with mean $0$ and variance $\sigma^2\geq 0$.  If $X$ follows this distribution, which we denote by $X\sim \cN_\C(0,\sigma^2)$, then $\Re X$ and $\Im X$ are independent  real Gaussian variables with  mean zero and variance $\sigma^2/2$. For $\sigma^2>0$, the Lebesgue density of $X$ is given by $(\pi\sigma^{2})^{-1}\eee^{- |z|^2/\sigma^2}$, $z\in\C$, whereas for $\sigma^2=0$ we have $X=0$ a.s.

\section{Main results}\label{sec:main_results}
Recall that $\xi_0,\xi_1,\ldots$ is a sequence of i.i.d.\ random variables taking complex values and that
we are interested in the  sequence of random polynomials
$$
Q_n(z) := \prod_{k=0}^n ( z - \xi_k).
$$
The next result describes the location of the critical points of $Q_n$ near its zero $\xi_0$.
We need the Cauchy--Stieltjes transform of the $\xi_k$'s, which is defined by
\begin{align}
f(z)=
\E \left[\frac{1}{z-\xi_1}\right] \label{eq:def.Cauchytr2.}
\end{align}
for those values $z\in\C$ for which $\E |\frac{1}{z-\xi_1}| <\infty$.
\begin{theorem} \label{theo:criticalpoints_random_u0}
Assume that the $\xi_k$'s have a Lebesgue density $p \colon \C \to [0,\infty)$ which is continuous on some open set $\mathscr D\subset \C$ and vanishes on $\C\backslash \mathscr D$.  Further, let $f$ be non-zero Lebesgue-a.e.\ on $\mathscr D$. Finally, let $r_1,r_2,\ldots >0$ be a positive sequence satisfying
\begin{align}
\lim_{n\to\infty} n r_n = + \infty
\quad
\text{and}
\quad
\lim_{n\to\infty} \sqrt{n} r_n = 0.
\end{align}
\begin{itemize}
\item[(a)]
The probability of the event that $Q_n$ has exactly one critical point in the disk  $\bB_{r_n}(\xi_0)$ converges to $1$ as $n\to\infty$, namely
$$
\lim_{n \to \infty} \P \left[\text{there is unique }  \zeta \in \bB_{r_n}(\xi_0) \text{ such that } Q_n'(\zeta)=0  \right] =1.
$$
\item[(b)]
Denoting the unique critical point of $Q_n$ in $\bB_{r_n}(\xi_0)$ by $\zeta_n$, if it exists uniquely, and defining $\zeta_n$ arbitrarily otherwise, we have
\begin{align}
\frac{f^2 (\xi_0)}{\sqrt{\pi p(\xi_0)}} \sqrt{\frac{n}{\log n}}   \left( n \left( \zeta_n-\xi_0 \right) +\frac{1}{f \left( \xi_0 \right)}\right)
\todistr \mathcal N_\C(0,1).
\end{align}
\end{itemize}
\end{theorem}
As we shall show in Lemma~\ref{lem:lipschitz}, below, the Cauchy--Stieltjes transform $f$ exists finitely everywhere on $\mathscr D$.
The above theorem will be deduced from the following somewhat easier statement in which $\xi_0$ is replaced by a deterministic zero. Let $\xi_1,\xi_2,\ldots$ be i.i.d.\ random variables with complex values. Fix some deterministic $u_0\in \C$ and consider the random polynomials
$$
P_n(z) := (z-u_0) \prod_{k=1}^n ( z - \xi_k ).
$$
We are interested in the location of the critical point near $u_0$.
\begin{theorem} \label{theo:criticalpoints}
Assume that on a sufficiently small disk around $u_0$, the random variables $\xi_k$ have a Lebesgue density $p$ that is continuous at $u_0$. Also, let the Cauchy--Stieltjes transform $f$ given by~\eqref{eq:def.Cauchytr2.}  satisfy $f(u_0) \neq 0$. Finally, let  $r_1,r_2,\ldots > 0$ be a sequence of positive numbers satisfying
\begin{align}
\lim_{n\to\infty} n r_n = + \infty
\quad
\text{and}
\quad
\lim_{n\to\infty} \sqrt{n} r_n = 0.
\end{align}
\begin{itemize}
\item[(a)]
The probability of the event that the random polynomial $P_n$ has exactly one critical point in the disk $\bB_{r_n}(u_0)$ converges to $1$ as $n \to \infty$, that is
\begin{align}
\lim_{n \to \infty} \P [\text{there is unique }  \zeta \in \bB_{r_n}(u_0) \text{ such that }  P_n'(\zeta)=0] = 1. \label{eq:mainthm1}
\end{align}
\item[(b)]
Denoting this critical point by $\zeta_n$, if it exists uniquely, and defining $\zeta_n$ arbitrarily otherwise, we have
\begin{align}
f^2 ( u_0 ) \sqrt{\frac{n}{\log n}}
\left( n \left( \zeta_n-u_0 \right) +\frac{1}{f \left( u_0 \right)}\right)  \todistr \cN_{\C} \left( 0 , \pi p(u_0) \right) \text{.} \label{eq:mainthm2}
\end{align}
\end{itemize}
\end{theorem}
Note that outside a small neighborhood of $u_0$ the distribution of $\xi_1$ may be completely arbitrary, for example it may have atomic or singular components.  The existence of the Cauchy--Stieltjes transform at $u_0$ will be established in Lemma~\ref{lem:lipschitz}, below.

It follows from Theorem~\ref{theo:criticalpoints}~(a) that for every $\alpha \in (1/2,1)$, the probability that the disk $\bB_{n^{-\alpha}}(u_0)$ contains exactly one critical point $\zeta_n$ converges to $1$. On the other hand, one can prove that the distance between $u_0$ and the closest zero of $P_n$ satisfies
\begin{equation}\label{eq:weibull}
\lim_{n\to\infty} \P\left[\sqrt n \min_{k=1,\ldots,n} |\xi_k - u_0| \leq r \right] = 1- \eee^{-p(u_0)\pi r^2},
\quad
r>0.
\end{equation}
Indeed, since the  density of the $\xi_k$'s is continuous at $u_0$, we have $\P[ |\xi_k - u_0| \leq s] \sim p(u_0) \pi s^2$ as $s\downarrow 0$, from which~\eqref{eq:weibull}  easily follows.
That is,  in the case when $p(u_0) >0$, the typical distance from $u_0$ to the closest zero is of order $1/\sqrt n$. As we shall see in the next paragraph,  the distance to the associated critical point $\zeta_n$ is of order $1/n$, which is much smaller.

Part (b) of the theorem describes the ``fluctuations'' of the critical point $\zeta_n$. Roughly speaking, part (b) states that $\zeta_n$ satisfies
\begin{equation}\label{eq:zeta_n_approx}
\zeta_n = u_0 - \frac {1}{n f(u_0)}  + \frac{\sqrt{\pi p(u_0)\log n}}{n^{3/2} |f^2(u_0)|} (N+ o(1)),
\end{equation} 
where $N\sim \cN_\C(0,1)$ is complex standard normal. The next corollary of Theorem~\ref{theo:criticalpoints} provides a ``confidence disk'' for the critical point associated with $u_0$.

\begin{corollary}
For every fixed $R>0$, the probability that the critical point $\zeta_n$ is not contained in the disk of radius
$$
\frac{\sqrt{\pi p(u_0)\log n}}{n^{3/2} |f^2(u_0)|} R
$$
centered at the point $u_0 - \frac 1 {n f(u_0)}$ converges to $\eee^{-R^2}$, as $n\to\infty$.
\end{corollary}
\begin{proof}
The probability mentioned in the statement of the corollary equals
$$
\P\left[  \sqrt{\frac n {\log n}} \frac{|f^2 (u_0)|}{\sqrt{\pi p(u_0)}} \left| n \left( \zeta_n-u_0 \right) +\frac{1}{f \left(u_0 \right)}\right|  > R\right].
$$
By Theorem~\ref{theo:criticalpoints} (b) and the continuous mapping theorem, this expression converges to $\P[|N| > R] = \eee^{-R^2}$, where $N\sim \cN_\C(0,1)$ is complex standard normal random variable.
\end{proof}

\section{Examples}\label{sec:examples}
In this section we shall give several special cases of the above results. In these examples, the density $p(z)$ is rotationally invariant, and the Cauchy--Stieltjes transform defined in~\eqref{eq:def.Cauchytr2.} can be computed explicitly. The next proposition, which is standard, follows essentially from the fact that the two-dimensional electrostatic field generated by the uniform probability distribution on a circle centered at $0$ vanishes inside the circle and coincides with the field generated by a unit charge at $0$ outside the circle.
\begin{proposition} \label{prop:Cauchy p(abs)}
Let $\xi_1$ be a complex random variable whose Lebesgue density $p(z)=q(|z|)$ depends on its argument $z$ only by its absolute value.
Then the Cauchy--Stieltjes transform of $\xi_1$ is
\begin{align}
f(z)= \frac{1}{z} \P \left[|\xi_1| \le |z| \right]
= \frac{2 \pi}{z} \int_0^{|z|} r q(r) \dd r,
\quad
\text{ if } z \neq 0,
\end{align}
and $f(0)=0$.
\end{proposition}
\begin{proof}
Let first  $z \neq 0$.
Using the identity
\begin{align}
\int_0^{2 \pi} \frac{\dd \phi}{z-r\eee^{i \phi}}=
\begin{cases}
0, & \text{ if } |z|<r, \\
\frac{2 \pi}{z}, & \text{ if } |z|>r,
\end{cases} 
\end{align}
after  passing to polar coordinates, we obtain
\begin{multline*}
f(z)=\int_\C \frac{p(u)}{z-u} \dd u
=\int_0^\infty r q(r) \int_0^{2 \pi} \frac{1}{z-r\eee^{i \phi}} \dd \phi \dd r
=\frac{2 \pi}{z} \int_0^\infty r q(r) \indi{|z|>r} \dd r
\\
=\frac{2 \pi}{z} \int_0^{|z|} r q(r) \dd r
=\frac{1}{z} \int_{\bB_{|z|}(0)} p(u) \dd u
=\frac{1}{z} \P \left[|\xi_1| \le |z|\right].
\end{multline*}
For $z=0$ we have, using polar coordinates,
\begin{align*}
f(0)=- \int_{\C} \frac{p(u)}{u} \dd u
=-\int_0^\infty r q(r) \int_0^{2 \pi} \frac{1}{r \eee^{i \phi}} \dd \phi \dd r =0 \text{.}
\end{align*}
\end{proof}
\begin{example}\label{eq:uniform_iid_zeroes}
If $\xi_0,\xi_1,\ldots$ are i.i.d.\ with the uniform distribution on the unit disk, the density is given by  $p(z) := \frac{1}{\pi} \indi{|z|<1}$ and the Cauchy--Stieltjes transform is
\begin{align*}
f(z)
=\frac{1}{z} \P \left[|\xi_1| \le |z|\right]
=\frac{1}{z} \left( \min \left\{ |z|,1 \right\} \right)^2
=
\begin{cases}
\overline{z}, &\text{ if } |z|\leq 1,\\
1/z, &\text{ if } |z|\geq 1.
\end{cases}
\end{align*}
Since $|\xi_0|<1$ a.s.,\ Theorem~\ref{theo:criticalpoints_random_u0} takes the form
$$
\bar \xi_0^2 \sqrt {\frac{n}{\log n}}  \left(n (\zeta_n - \xi_0) + \frac{1}{\bar \xi_0}\right) \todistr \Normal_\C (0,1).
$$
In fact, $\bar \xi_0^2$ can be replaced with $|\xi_0|^2$ since the result holds conditionally on $\xi_0$ (Theorem~\ref{theo:criticalpoints}) and $\Normal_\C (0,1)$ is a rotationally invariant distribution.
Note that $1/\bar \xi_0$ becomes large if $\xi_0$ is close to $0$, which explains why the distance between the zero and the corresponding critical point tends to become larger for zeroes close to the origin; see Figure~\ref{fig:zeroes_critical}.
\end{example}

\begin{example}
If $\xi_1,\xi_2,\ldots$ are standard complex normal, the density is $p(z)=\frac{1}{\pi} \eee^{-|z|^2}$ and the Cauchy--Stieltjes transform is
\begin{align*}
f(z) = \frac{2}{z} \int_0^{|z|} r \eee^{-r^2} \dd r
=\frac{1}{z} \left( 1-\eee^{-|z|^2} \right) \text{,}
\end{align*}
if $z \neq 0$, and $f(0) = 0$. 
\end{example}

\begin{figure}[!tbp]
\centering
\includegraphics[width=0.49\textwidth]{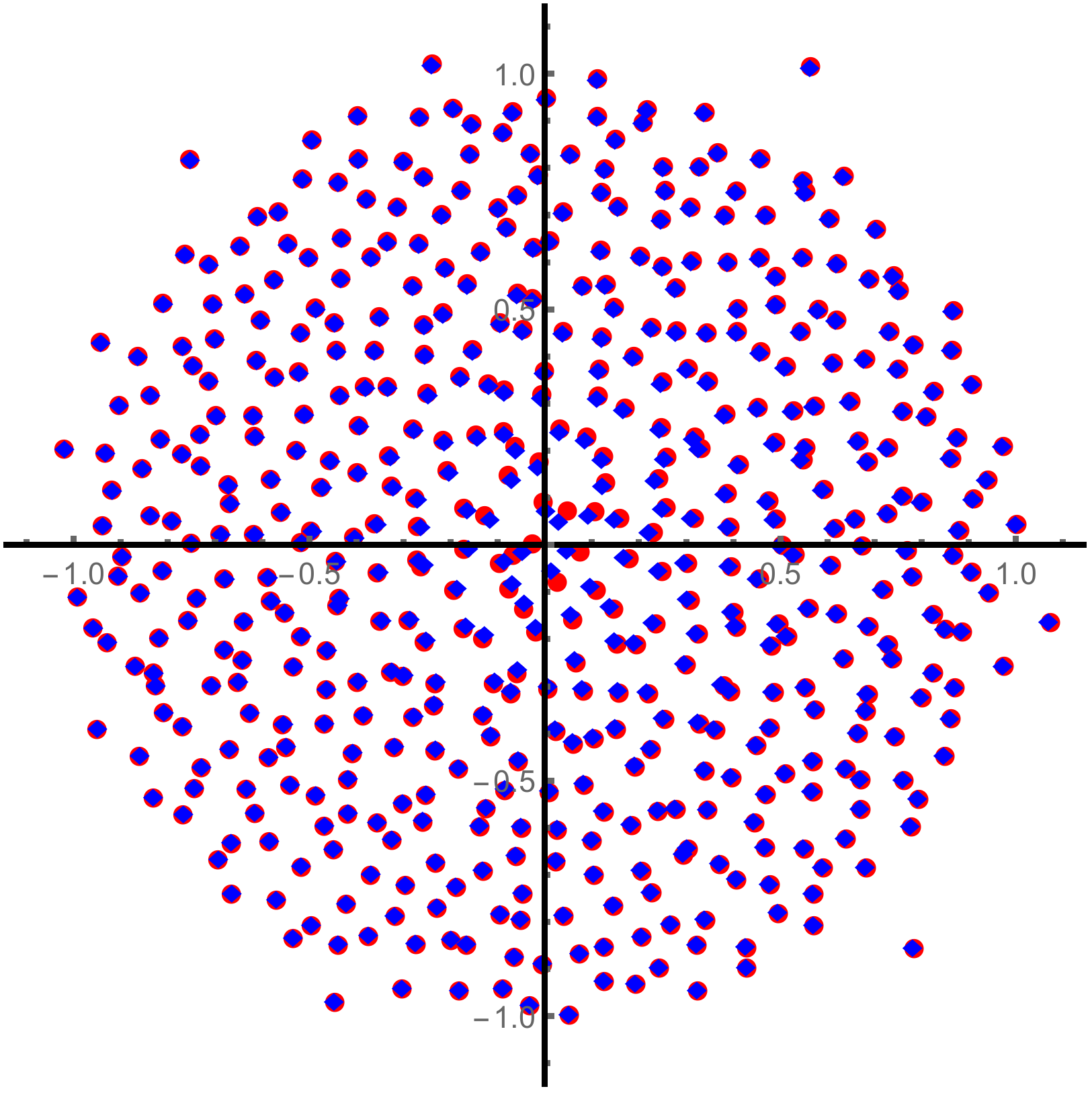}
\includegraphics[width=0.49\textwidth]{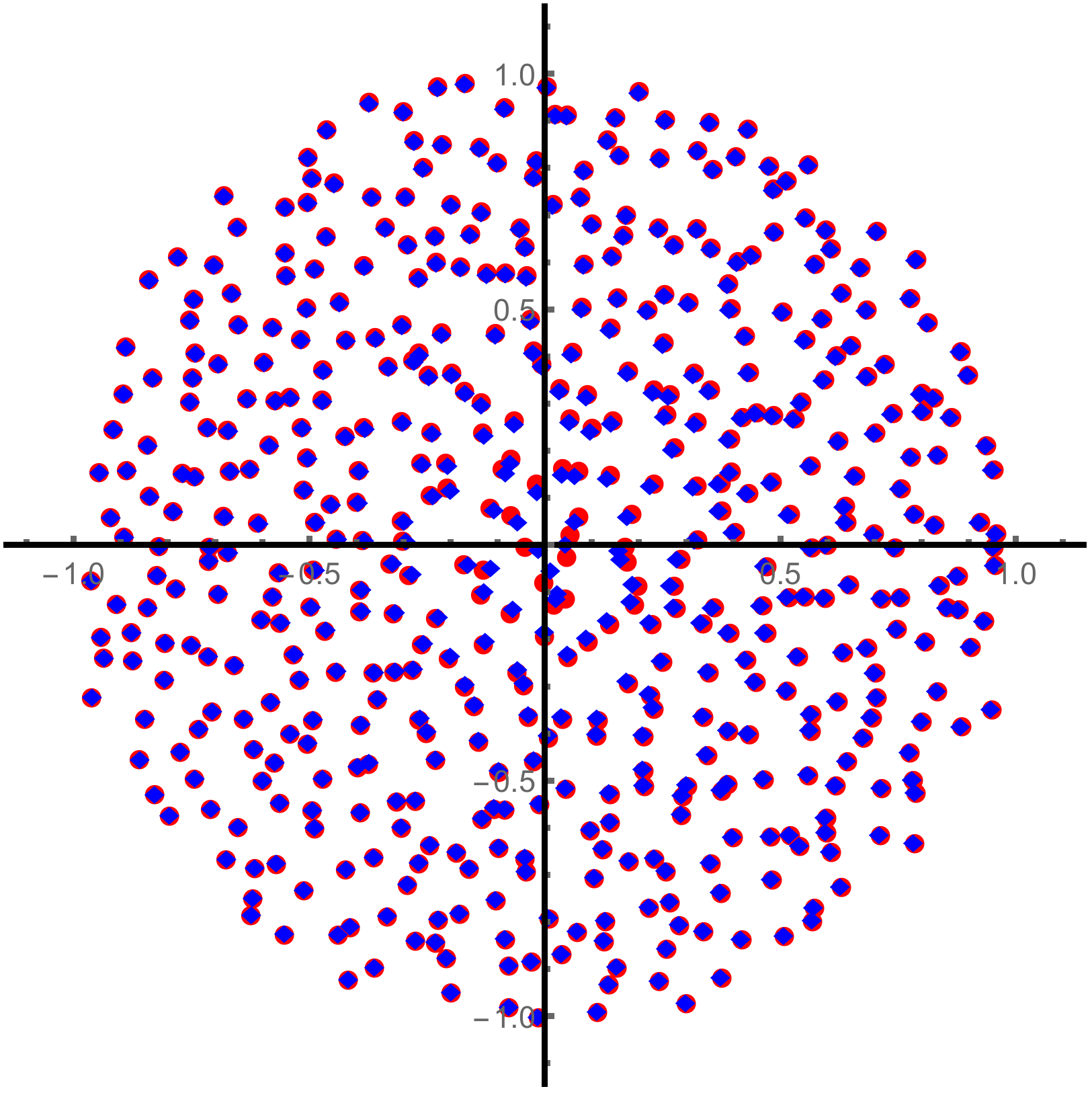}
\caption{Zeroes and critical points of a Weyl polynomial (left) and the characteristic polynomial of a Ginibre random matrix (right). The degree is $n=500$ in both cases.  Red disks: zeroes. Blue diamonds: critical points.}
	\label{fig:zeroes_critical_weyl_ginibre}
\end{figure}

\section{Conjectures}\label{sec:conjectures}
The above results suggest several conjectures on critical points of random polynomials with not necessarily i.i.d.\ zeroes.
Roughly speaking, we shall try to transform the information on the quality of the pairing near some individual zero (see Theorem~\ref{theo:criticalpoints_random_u0}) into conjectures on the behavior of the distances between all zeroes and the corresponding critical points.
For concreteness, we shall consider the following three families of random polynomials whose zeroes are asymptotically uniformly distributed  on the unit disk.
\begin{enumerate}
\item[(a)]
Polynomials with i.i.d.\ zeroes
$$
p_n^{\text{i.i.d.}}(z) = \prod_{k=1}^n (z - \xi_k),
$$
where $\xi_1,\xi_2,\ldots$ are i.i.d.\ random variables having uniform distribution on the unit disk; see Figure~\ref{fig:zeroes_critical}.
\item[(b)] Weyl polynomials
$$
p_n^{\text{Weyl}} (z) = \sum_{k=0}^n X_k \frac{(z \sqrt n)^k}{\sqrt {k!}},
$$
where $X_0,X_1,\ldots$ are i.i.d.\ random variables; see Figure~\ref{fig:zeroes_critical_weyl_ginibre}, left panel.
\item[(c)]
Characteristic polynomials of the form
$$
p_n^{\text{char}} (z) = \det (A_n - z \sqrt n),
$$
where $A_n = (a_{ij})_{i,j=1}^n$ is a random $n\times n$-matrix with i.i.d.\ entries $a_{ij}$; see Figure~\ref{fig:zeroes_critical_weyl_ginibre}, right panel.
\end{enumerate}
Let $p_n(z)$ be a random polynomial chosen according to one of the above models. Denote by  $Z_{1,n},\ldots, Z_{n,n}$ its complex zeroes, and let $W_{1,n}, \ldots, W_{n-1,n}$ be its critical points, i.e.\ the zeroes of $p_n'$. The empirical distribution of zeroes and the empirical distribution of critical points are random probability measures on $\C$ defined by
$$
\mu_n^{\text{zeroes}} := \frac 1n \sum_{k=1}^n \delta_{Z_{k,n}},
\qquad
\mu_n^{\text{crit}} := \frac 1{n-1} \sum_{k=1}^{n-1} \delta_{W_{k,n}},
$$

It is known that for all three models, the empirical distribution of zeroes
converges to the uniform distribution on the unit disk provided suitable moment conditions are satisfied. More precisely, let $\mathcal M_1(\C)$ be the space of probability measures on $\C$ endowed with the weak topology. We say that a sequence $(\mu_n)_{n\in\N}$ of random elements with values in $\mathcal M_1(\C)$ converges in probability to some deterministic probability measure $\mu\in \mathcal M_1(\C)$, if for every $\eps>0$,
$$
\lim_{n\to\infty} \P[ \rho(\mu_n,\mu) > \eps] = 0,
$$
where $\rho$ is any metric generating the weak topology on $\mathcal M_1(\C)$. This mode of convergence is denoted by $\mu_n\toprobabnon \mu$.
For all three models (under appropriate moment conditions) we have
$$
\mu_n^{\text{zeroes}} \toprobabnon \text{Unif}(\bB_1(0)), \qquad \text {as } n\to\infty,
$$
where $\text{Unif}(\bB_1(0))$ is the uniform probability distribution on the unit disk. For polynomials with i.i.d.\ zeroes, this is just the law of large numbers for empirical processes. For characteristic polynomials, this is the circular law; see~\cite{tao_vu} for the the proof under the assumption $\E a_{ij}=0$, $\E |a_{ij}|^2 = 1$. Finally, for Weyl polynomials, this was proved in~\cite{kabluchko_zaporozhets12a} under the assumption $\E \log_+ |X_0|<\infty$.
Similarly, for the empirical measure of critical points it is known that
$$
\mu_n^{\text{crit}} \toprobabnon \text{Unif}(\bB_1(0)), \qquad \text {as } n\to\infty,
$$
in the case of Weyl polynomials~\cite[Remark~2.11]{kabluchko_zaporozhets12a} and polynomials with i.i.d.\ zeroes~\cite{kabluchko15}. For characteristic polynomials, this was conjectured in~\cite{orourke_matrices}, where analogous relations were established for some other ensembles of random matrices.

\begin{figure}[!tbp]
\centering
\includegraphics[width=0.32\textwidth]{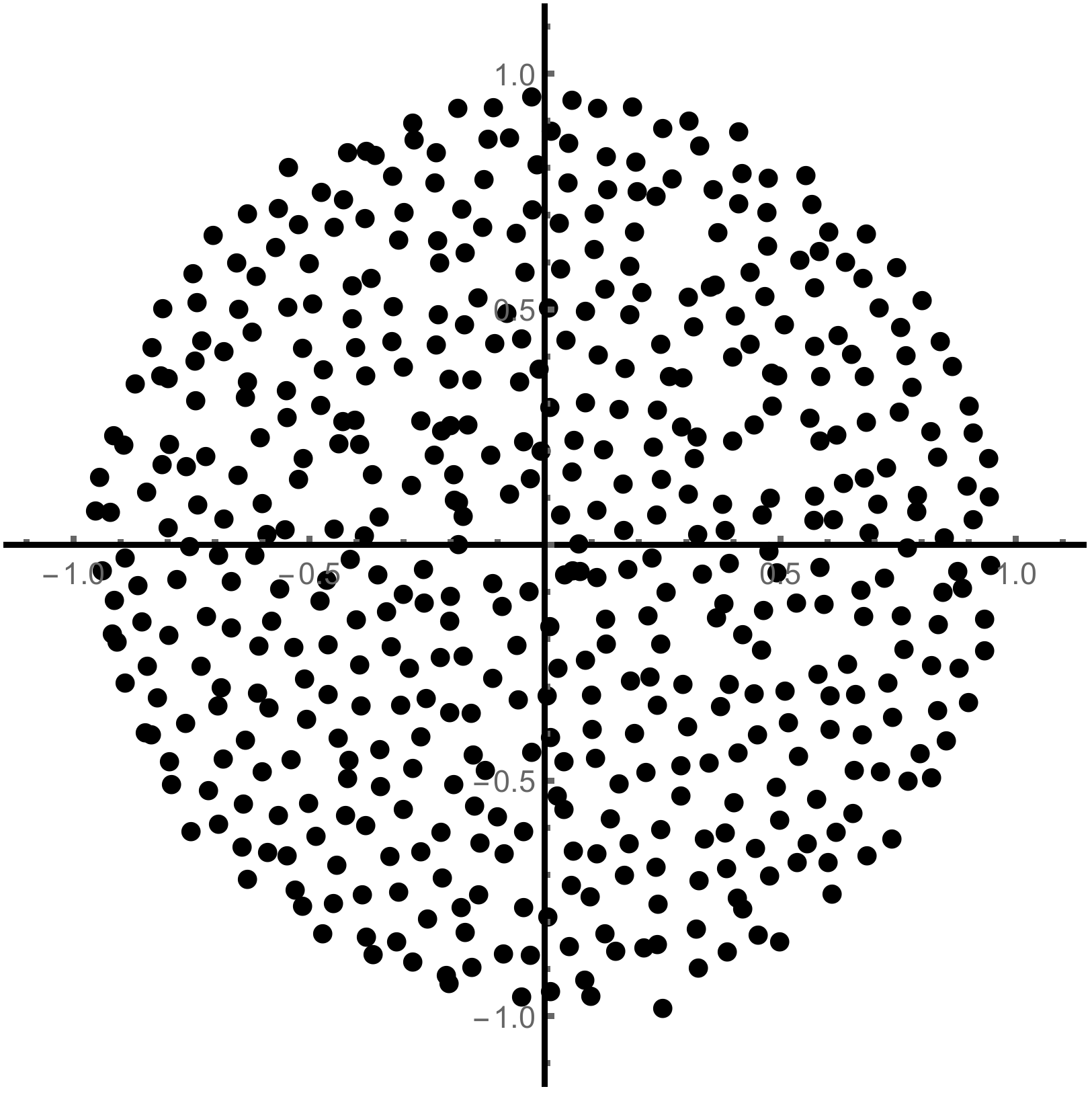}
\includegraphics[width=0.32\textwidth]{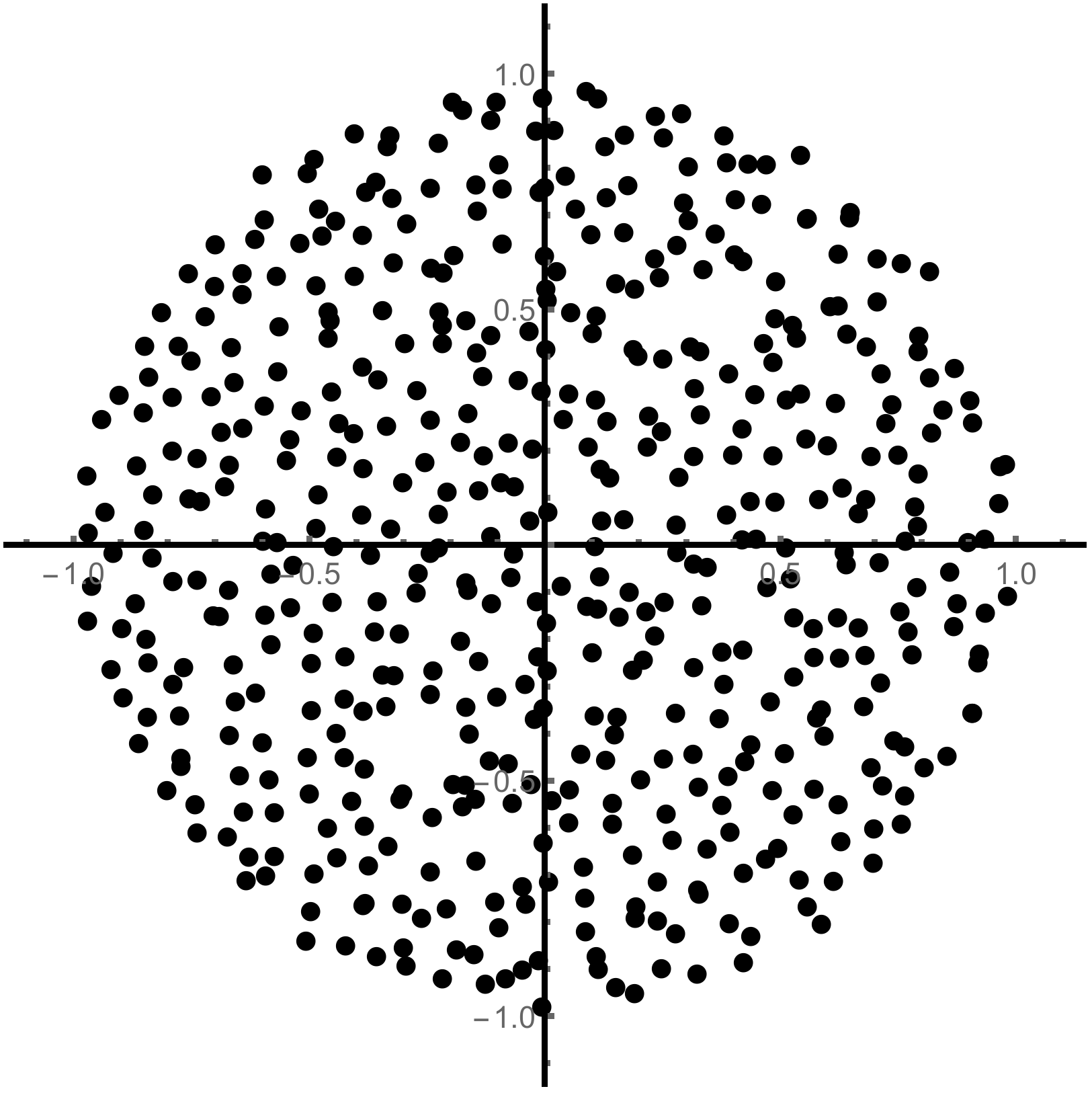}
\includegraphics[width=0.32\textwidth]{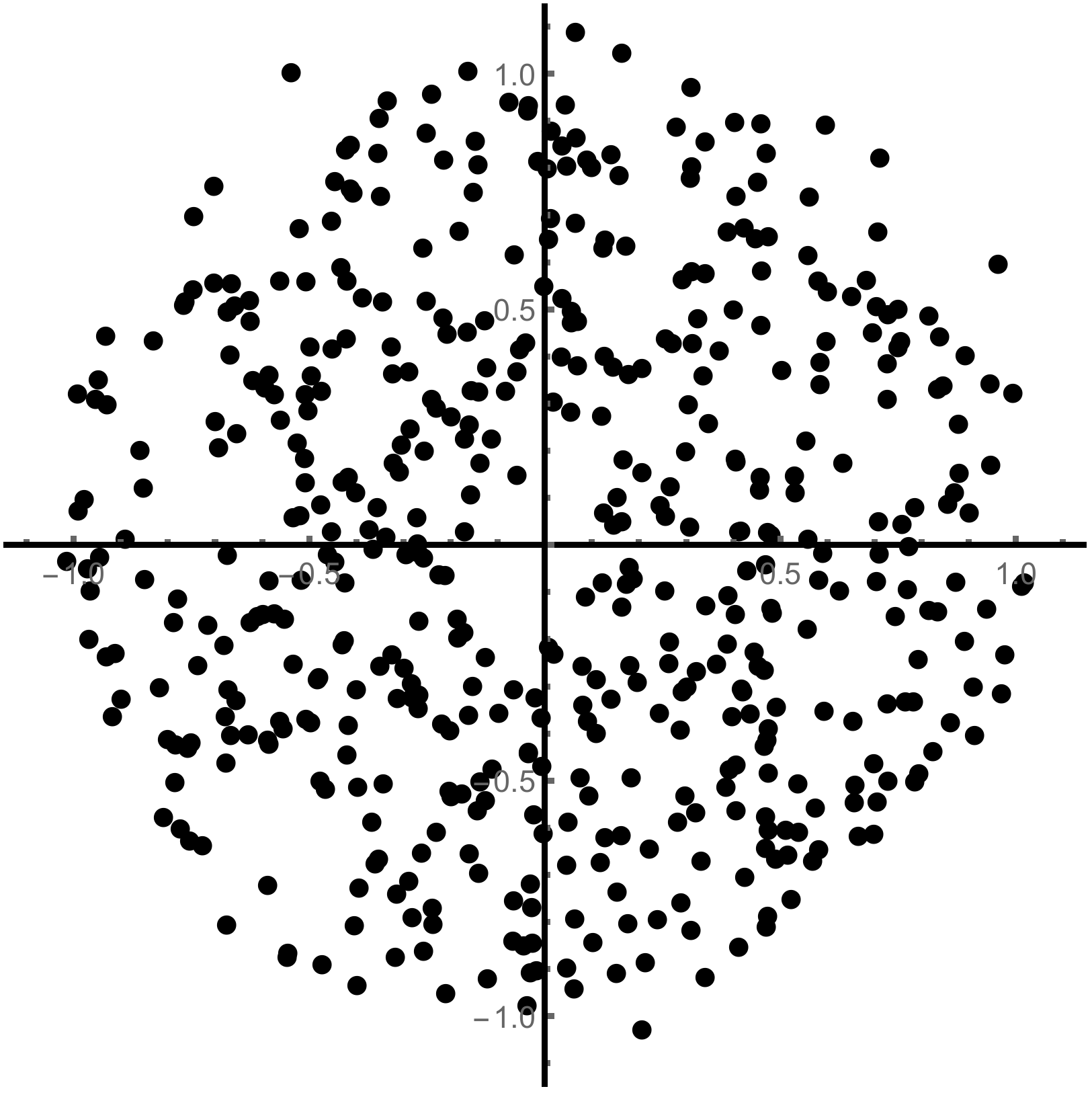}
\caption{The sample $\{ 1 / (n (Z_{k,n} - \zeta_{k,n}))\colon k=1,\ldots,n\}$ of the normalized inverse differences between the zeroes and the associated critical points. Left panel: Weyl polynomial. Middle panel: Characteristic polynomial of a Ginibre random matrix. Right panel: Polynomial with i.i.d.\ zeroes distributed uniformly on the unit disk.  The degree is $n=500$ in all three cases.}
	\label{fig:distances_rough_scale}
\end{figure}

\begin{figure}[!tbp]
\centering
\includegraphics[width=0.35\textwidth]{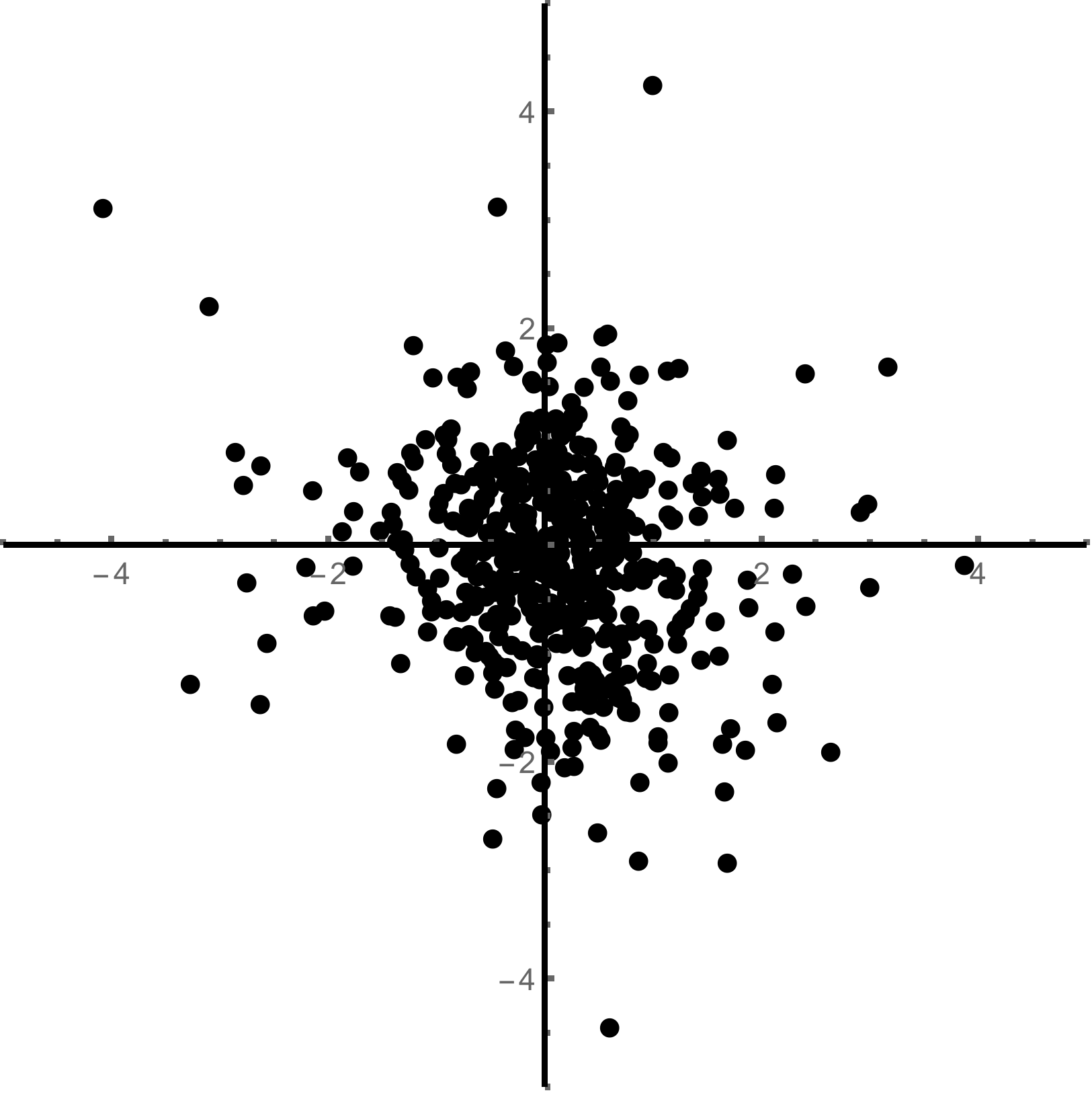}
\includegraphics[width=0.55\textwidth]{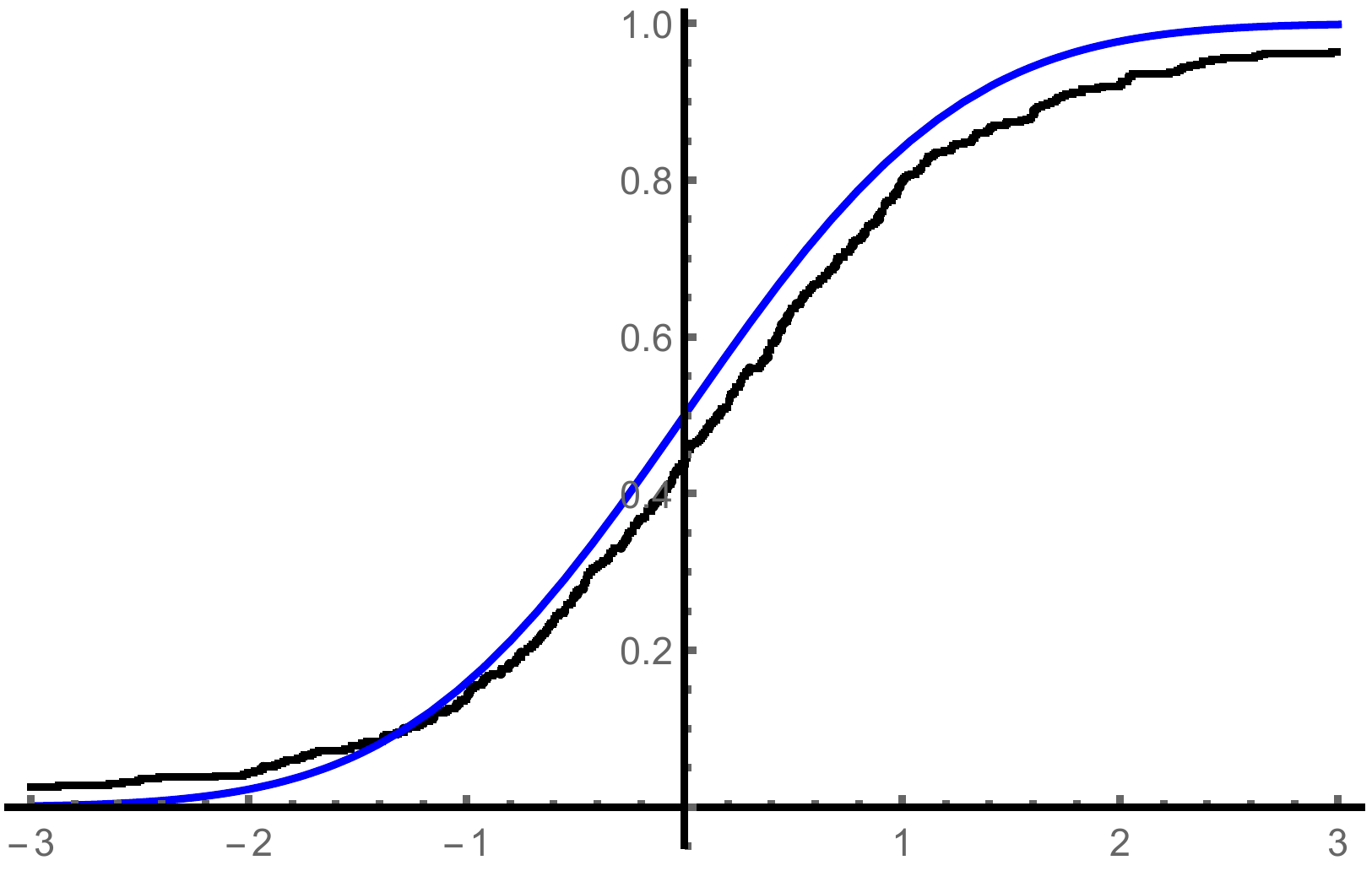}
\caption{Left: Atoms in a realization of $\nu_n^{(1)}$ for $Q_n^{\text{i.i.d.}}$. Right: Empirical distribution function (black) of the sample $\{\sqrt 2 \Re d_{n,k}\colon 1\leq k\leq n\}$ for the  random polynomial $p_n^{\text{i.i.d.}}$, together with the standard normal distribution function (blue). The degree is $n=500$ in both cases.}
	\label{fig:distances}
\end{figure}

To state conjectures on the quality of the pairing between the zeroes and the critical points, we need to introduce more notation. For every zero $Z_{k,n}$, where $k\in\{1,\ldots,n\}$, let $\zeta_{k,n}$ be the critical point of $p_n$ most close to $Z_{k,n}$, that is $p_n'(\zeta_{k,n}) = 0$ and
$$
|\zeta_{k,n} - Z_{k,n}| = \min \{|z-Z_{k,n}|\colon z\in\C, p_n'(z)=0\}.
$$
For polynomials with i.i.d.\ zeroes, Theorem~\ref{theo:criticalpoints} (see also~\eqref{eq:zeta_n_approx} and Example~\ref{eq:uniform_iid_zeroes}) suggests the approximations
\begin{equation}\label{eq:approximation}
\zeta_{k,n} \approx Z_{k,n} - \frac {1} {n \bar Z_{k,n}},
\quad
\zeta_{k,n} \approx Z_{k,n} - \frac {1}{n \bar Z_{k,n}}  + \frac{\sqrt{\log n}}{n^{3/2}|Z_{k,n}|^2} N_k,
\end{equation}
where $N_k\sim \mathcal N_\C(0,1)$.  In order to quantify the quality of these approximations, introduce the random measures
$$
\nu_n
:=
\sum_{k=1}^n \delta_{1/(n (Z_{k,n} - \zeta_{k,n}))}
\quad
\text{and}
\quad
\chi_n
:=
\sum_{k=1}^n \delta_{d_{k,n}},
$$
where
\begin{equation}\label{eq:def_d_k_n}
d_{k,n} := |Z_{k,n}|^2 \sqrt{\frac n {\log n}} \left(n(Z_{k,n} - \zeta_{k,n}) - \frac 1 {\bar Z_{k,n}}\right) .
\end{equation}
We view $\nu_n$ and $\chi_n$ as random elements with values in the space $\mathcal M_1(\C)$; see Figure~\ref{fig:distances_rough_scale} for realisations of $\nu_n$ in all three models and Figure~\ref{fig:distances} (left panel) for a realisation of $\chi_n$ in the i.i.d.\ zeroes model. If the first approximation in~\eqref{eq:approximation} is valid, then we should have $1/(n (Z_{k,n} - \zeta_{k,n})) \approx \bar Z_{k,n}$. This suggests that $\nu_n$ should be close to the distribution of the $\bar Z_{k,n}$'s, which is uniform on the unit disk. Numerical simulations, see Figure~\ref{fig:distances_rough_scale}, support the following

\begin{conjecture}\label{conj:1}
For all three models, under appropriate moment conditions, we have
$$
\nu_n \toprobabnon \text{Unif}(\bB_1(0)), \qquad \text {as } n\to\infty.
$$
\end{conjecture}

Similarly, if the second, more refined approximation in~\eqref{eq:approximation} is valid, then we should have $d_{k,n}\approx -N_k$. The following conjecture is supported by numerical simulations:
\begin{conjecture}\label{conj:2}
For random polynomials with i.i.d.\ zeroes  we have
$$
\chi_n \toprobabnon \mathcal N_{\C} (0,1), \qquad \text {as } n\to\infty,
$$
where $\mathcal N_{\C} (0,1)$ denotes the standard normal distribution on $\C$.
\end{conjecture}

Numerical simulations suggest that the analogue of Conjecture~\ref{conj:2} holds for characteristic polynomials if $\sqrt {n /\log n}$ is replaced by another (unknown) normalizing sequence in the definition of $d_{k,n}$; see~\eqref{eq:def_d_k_n}. For Weyl polynomials, the limit seems to be some heavy-tailed, non-normal distribution.


\section{Proofs}\label{sec:proofs}
If no cancellation occurs, the critical points of $P_n$ coincide with the zeroes of its logarithmic derivative
$$
\frac{P_n'(z)}{P_n(z)} = \frac{1}{z-u_0} +\sum_{k=1}^n \frac{1}{z-\xi_k}.
$$
To prove Theorem~\ref{theo:criticalpoints}, we shall establish a functional limit theorem for the logarithmic derivative in a suitable small scaling window  near the conjectured location of the critical point. Then, we shall show  that the zero of the limit process (which is easy to compute) approximates the rescaled  critical point of $P_n$. Since the logarithmic derivative is a sum of independent random variables, it is natural to conjecture that its functional limit is a Gaussian process. This is indeed true. However, it turns out that the second moment of the summands is (just) infinite: they are in the non-normal domain of attraction of the normal law. This explains the logarithmic factor appearing in Theorems~\ref{theo:criticalpoints_random_u0} and~\ref{theo:criticalpoints}.

\subsection{Functional central limit theorem}
The classical central limit theorem states that if $X_1, X_2,\ldots$ are non-degenerate i.i.d.\ real-valued random variables with finite second moment, then
\begin{equation}\label{eq:CLT_classical}
\frac{X_1 + \ldots + X_n - a_n}{b_n} \todistr \cN(0,1),
\end{equation}
where we may take $a_n = n\E X_1$ and $b_n = \sqrt{n \Var X_1}$. However, there exist examples of random variables with \textit{infinite} second moment for which~\eqref{eq:CLT_classical} continues to hold with the same limiting normal distribution but with some other choice of the normalizing sequences $a_n$ and $b_n$. The set of all such variables, referred to as the \textit{non-normal domain of attraction} of the normal distribution, can be characterized by the condition $\P[|X_1|>t] = \ell(t) t^{-2}$ for some function $\ell(t)$ that varies slowly at $+\infty$; see~\cite[Theorem~2.6.2]{ibragimov_linnik_book}.

We shall be interested in the complex version of this situation. The following example is one of the simplest ones. If $\xi_1,\xi_2,\ldots$ are independent random variables having uniform distribution on the unit disk $\bB_1(0)$ in the complex plane, then we claim that
\begin{equation}\label{eq:example_CLT}
\frac {1}{\sqrt{n\log n}} \sum_{k=1}^n \frac 1 {\xi_k} \todistr \cN_{\C} (0, 1).
\end{equation}
Observe that the second moment of $\frac 1 {\xi_1}$ is infinite, and notice the additional factor $\sqrt {\log n}$ in the normalization. This term has the same origin as the factor which showed up in Theorems~\ref{theo:criticalpoints_random_u0} and~\ref{theo:criticalpoints}.
In fact, we shall need a more general theorem.


\begin{theorem} \label{theo:CLT}
Let $\xi_1,\xi_2,\ldots$ be a sequence of i.i.d.\ random variables with complex values. Assume that in some neighborhood of $0$,  these variables have Lebesgue density $p$ that is continuous at $0$.
Let $(z_n)_{n \in \N}$ be a complex sequence satisfying  $|z_n| = O(n^{-1/2 - \kappa})$, as $n\to\infty$, for some $\kappa>0$.
Then,
\begin{align}
\frac{1}{\sqrt{n \log n}} \sum_{k=1}^n \left( \frac{1}{z_n-\xi_k} - f(0) \right) \todistr \cN_\C \left(0,
\pi p(0)
\right).
\end{align}
\end{theorem}
For example, if the $\xi_k$'s have the  uniform distribution on the unit disk, we arrive at~\eqref{eq:example_CLT} by taking $z_n=0$ and observing that $f(0)=0$. The sequence $z_n$ in Theorem~\ref{theo:CLT} was introduced for technical reasons which will become clear later.

Next we state a functional version of Theorem~\ref{theo:CLT}. It will play a crucial role in our proof of Theorem~\ref{theo:criticalpoints}. For $\rho>0$ we denote by $\bA_\rho$ the space of continuous functions on the closed disk $\bB_\rho(0) = \{|z|\leq \rho\}$ which are analytic on the open disk $\{|z| < \rho\}$. Endowed with the supremum norm, $\bA_\rho$ becomes a Banach space.
\begin{theorem} \label{theo:funct.CLT}
Let $\xi_1,\xi_2,\ldots$ be a sequence of i.i.d.\ random variables with complex values. Assume that in some neighborhood of $0$,  these variables have Lebesgue density $p$ that is continuous at $0$, and let their Cauchy--Stieltjes transform $f$ satisfy $f(0)\neq 0$.
Fix any $\rho >0$ and define the functions
\begin{align}\label{eq:z_n_w_def}
z_n(w) := - \frac{1}{n f(0)} + \frac{w \sqrt{\log n}}{n^\frac{3}{2} (f(0))^2 }, \quad |w| \leq \rho.
\end{align}
Then we have the following weak convergence of stochastic processes on the space $\bA_\rho$:
\begin{align}\label{theo:funct.CLT.conv}
\left(\frac{1}{\sqrt{n \log n}} \sum_{k=1}^n \left( \frac{1}{z_n (w) -\xi_k} - f(0) \right)\right)_{|w|\leq \rho}
\toweak
(N)_{|w|\leq \rho},
\end{align}
where $N\sim \cN_\C (0, \pi p(0))$.
\end{theorem}
Note that each particular realisation of the limit process in~\eqref{theo:funct.CLT.conv} is a constant function.

\begin{remark}\label{rem:convention}
Observe also that the function on the left-hand side of~\eqref{theo:funct.CLT.conv} may have poles in the disk $\bB_\rho(0)$, in which case it is not an element of the space $\bA_\rho$. To turn the left-hand side of~\eqref{theo:funct.CLT.conv} into a well-defined random element of $\bA_\rho$, we agree to re-define the function to be identically $0$ each time it has a pole in $\bB_\rho(0)$.  Thus, instead of a meromorphic function $h$  we in fact consider the analytic function $h \indi{h\in\bA_\rho}$.  We shall see below that the probability of having poles in $\bB_{\rho}(0)$ goes to $0$ as $n\to\infty$, so that these poles do not affect the distributional convergence.
\end{remark}

\subsection{Proof of Theorem~\ref{theo:CLT}}

The proof is based on the following result, see~\cite[Theorem~3.2.2]{meerschaert_book}, which is
somewhat
more general than the Lyapunov (and even Lindeberg) central limit
theorem. We denote by $|\cdot|$
the Euclidean norm in $\R^m$ and by $\Cov Z$ the covariance matrix of an
$\R^m$-valued random vector $Z$.
\begin{theorem}[General CLT for Random
Vectors]\label{theo:rvaceva}\index{CLT!general}
For every $n\in\N$, let $\{\bZ_{n,k}\colon 1\leq k\leq k_n\}$ be independent
$\R^m$-valued random
vectors. Assume that the following conditions hold:
\begin{enumerate}

\item[\textup{(a)}] \label{cond:gned1} For every $\eps>0$,
$\lim_{n\to\infty} \sum_{k=1}^{k_n} \P[|\bZ_{n,k}|>\eps]=0$.

\item[\textup{(b)}] \label{cond:gned2} For some positive semidefinite
matrix $\Sigma$,
\begin{align}
\Sigma
=
\lim_{\eps\downarrow 0} \limsup_{n\to\infty} \sum_{k=1}^{k_n} \Cov
\left[\bZ_{n,k} \indi{|\bZ_{n,k}|<\eps}\right]
=
\lim_{\eps\downarrow 0} \liminf_{n\to\infty} \sum_{k=1}^{k_n} \Cov
\left[\bZ_{n,k} \indi{|\bZ_{n,k}|<\eps}\right].
\end{align}

\end{enumerate}
Then, the random vector $\bS_n:=\sum_{k=1}^{k_n} (\bZ_{n,k}-\E
[\bZ_{n,k}\ind_{\{|\bZ_{n,k}|<R\}}])$ converges weakly to a mean zero Gaussian
distribution on $\R^m$ with covariance matrix $\Sigma$. Here, $R>0$ is
arbitrary.
\end{theorem}

In the proof of Theorem~\ref{theo:CLT} we shall need the following
\begin{lemma} \label{lemma:Integralasymptote}
Fix some $R>0$ and let $g_{\delta} \colon [0,R] \to \C$ be a uniformly bounded family of complex-valued measurable functions parametrized by $\delta \in [0,1]$ and having the following property:
For every $\eps>0$ there exist $r_0, \delta_0 \in (0,R)$ such that
\begin{align*}
|g_\delta(r)-g_0(0)| <\eps \text{ for every } \delta \in [0,\delta_0] \text{ and } r\in [0,r_0].
\end{align*}
That is, $g_\delta(r)$ is continuous at $(0,0)$  as a function of $r$ and $\delta$.  Then, as $\delta \downarrow 0$, we have
\begin{align}
\lim_{\delta\downarrow 0} \frac{1}{\log \frac{1}{\delta}} \int_\delta^R \frac{g_\delta(r)}{r} \dd r
=
 g_0(0). \label{eq:lem.Integralasymptote-claim}
\end{align}
\end{lemma}

\begin{proof}
Fix any $\varepsilon >0$ and choose sufficiently small $\delta_0, r_0 \in (0,R)$ as in the lemma.
For every $\delta \in (0, \min\{ \delta_0, r_0,1 \})$ we have
\begin{align*}
\left|\int_\delta^{r_0} \frac{g_\delta(r)}{r} \dd  r
-
(\log r_0 -\log \delta)g_0(0)\right|
=
\left|\int_\delta^{r_0} \frac{g_\delta(r)-g_0(0)}{r} \dd  r
\right|
\leq
\int_\delta^{r_0} \frac{\eps}{r} \dd r
=
\eps (\log r_0-\log\delta).
\end{align*}
Dividing by $\log \frac 1 \delta > 0$, we obtain
\begin{align}
\left| \frac{1}{\log \frac{1}{\delta}}
\int_\delta^{r_0} \frac{g_\delta(r)}{r} \dd r
-
\frac{(\log r_0 -\log \delta)g_0(0)}{\log \frac{1}{\delta}}
 \right|
\le
\eps \frac{(\log r_0 -\log \delta)}{\log \frac{1}{\delta}}. \label{eq:lem asympt of integral main est}
\end{align}
Observing that
$
\frac{(\log r_0 -\log \delta)}{\log \frac{1}{\delta}}
$
converges to $1$ as $\delta\downarrow 0$,
we conclude that all accumulation points of $\frac{1}{\log \frac{1}{\delta}}
\int_\delta^{r_0} \frac{g_\delta(r)}{r} \dd r$, as $\delta\downarrow 0$, are contained in the disk of radius $\eps$ centered at $g_0(0)$.
Now,
$$
\limsup_{\delta\downarrow 0} \left|\int_{r_0}^R \frac{g_\delta(r)}{r} \dd r\right| < \infty
$$
by the assumption that the family $g_\delta$ is uniformly bounded. It follows that all accumulation points of the bounded function $\frac{1}{\log \frac{1}{\delta}}
\int_\delta^{R} \frac{g_\delta(r)}{r} \dd r$, as $\delta\downarrow 0$, are contained in the same disk of radius $\eps$ centered at $g_0(0)$. Since this holds for every $\eps>0$, we arrive at~\eqref{eq:lem.Integralasymptote-claim}.
\end{proof}

\begin{proof}[Proof of Theorem \ref{theo:CLT}]
We shall show that the random variables
$$
Z_{n,k} := \frac{1}{\left( z_n-\xi_k \right) \sqrt{n \log n}}, \quad  k=1,\ldots,n,
$$
fulfil the assumptions of Theorem~\ref{theo:rvaceva}. Observe that $Z_{n,1}, \ldots, Z_{n,n}$ are independent and have the same distribution as
$Z_n := Z_{n,1}$. These complex variables are considered as two-dimensional random vectors, so that $m=2$ in the setting of Theorem~\ref{theo:rvaceva}.

Whenever helpful we shall analyze $V_n:=z_n-\xi_1$ instead of $\xi_1$. We shall denote the density of $V_n$ by $q_n(z) =p(z_n - z)$.
Since $p$ is continuous at $0$, there is a sufficiently small $R>0$ such that $p$ exists and is bounded in the disk $\bB_{2R}(0)$ by some constant $C$.
Since $\limto{n}{\infty} z_n = 0$, the density  $q_n$ is bounded on the smaller disk $\bB_{R}(0)$ by the same constant $C$  provided $n$ is sufficiently large.

\medskip
\noindent
\textit{Condition (a) of Theorem \ref{theo:rvaceva}.} To verify this condition, we need to check that for every $\eta >0$,
\begin{equation}\label{eq:ver_Lindeberg_first_cond}
\limto{n}{\infty} \sum_{k=1}^{n} \P \left[ \left| \frac{1}{\left( z_n-\xi_k \right)\sqrt{n \log n} } \right| >\eta \right]
= \limto{n}{\infty} n \P \left[ \left| \frac{1}{V_n \sqrt{n \log n}} \right| >\eta \right] =0.
\end{equation}
For sufficiently large $n$, we have
\begin{equation*}
n  \P \left[ \left| \frac{1}{V_n \sqrt{n \log n}} \right| > \eta \right]
= n  \P \left[ \left| V_n \right| < \frac{1}{\eta \sqrt{n \log n}}\right]
= n \int_{\left\{ |z| < \frac{1}{\eta \sqrt{n \log n}} \right\} } q_n(z) \dd z
\le  \frac{n C \pi}{\eta^2 n \log n},
\end{equation*}
which converges  to $0$, thus verifying~\eqref{eq:ver_Lindeberg_first_cond}.

\medskip
\noindent
\textit{Condition (b) of Theorem \ref{theo:rvaceva}.}
To verify this condition, it suffices to show that for every $\eta>0$,
\begin{align}
\lim_{n \to \infty} \sum_{k=1}^n \Cov \left[ \frac{1}{\left( z_n-\xi_k \right) \sqrt{n \log n}} \indi{ \left| \frac{1}{\left( z_n-\xi_k \right) \sqrt{n \log n}} \right| < \eta} \right]
=
\begin{pmatrix}
\frac{\pi}{2} p(0) & 0 \\
0 & \frac{\pi}{2} p(0)
\end{pmatrix}. \label{eq:meerschaerd2_applied}
\end{align}
Recall that by identifying $\C$ with $\R^2$ we can consider complex-valued random variables as two-dimensional random vectors. The above can be simplified to
\begin{equation} \label{eq:meerschaerd2_applied_restate}
\lim_{n \to \infty} n \Cov \left[ Z_n \indi{ |Z_n| <\eta } \right]
=
\begin{pmatrix}
\frac{\pi}{2} p(0) & 0 \\
0 & \frac{\pi}{2} p(0)
\end{pmatrix}.
\end{equation}
It suffices  to prove that
\begin{align}
&\lim_{n\to\infty} n\E \left[ \left(\Re Z_n \indi{ \left| Z_n \right| < \eta } \right)^2\right]
=
\lim_{n\to\infty} n\E \left[ \left(\Im Z_n \indi{ \left| Z_n \right| < \eta } \right)^2\right]
=
\frac{\pi}{2} p(0), \label{eq:cov_lim1}
\\
&\lim_{n\to\infty}  n \E \left[ \Re Z_n \Im Z_n \indi{ \left| Z_n \right| < \eta } \right]
=0, \label{eq:cov_lim2}
\\
&\lim_{n\to\infty}  n (\E [\Re Z_n \indi{ \left| Z_n \right| < \eta }])^2
=
\lim_{n\to\infty}  n (\E [ \Im Z_n \indi{ \left| Z_n \right| < \eta }])^2
=0. \label{eq:cov_lim3}
\end{align}

\noindent
\textit{Proof of~\eqref{eq:cov_lim1} and~\eqref{eq:cov_lim2}.} Recall that $Z_n = 1/(V_n\sqrt{n\log n})$.
In view of the identities $\E Z \bar Z =\E \left( \Re Z \right)^2 + \E \left( \Im Z \right)^2$ and $\E Z^2 =\E \left( \Re Z \right)^2 - \E \left( \Im Z \right)^2 + 2 i \E \Re Z \Im Z$, it suffices to show that
\begin{equation}\label{eq:exp_z_n_overline_z_n}
\lim_{n\to\infty} n \E \left[Z_n \overline{Z_n}\indi{ |V_n| > \frac 1 {\eta \sqrt{n \log n}} }\right] = \pi p(0)
\quad
\text{and}
\quad
\lim_{n\to\infty} n \E \left[Z_n^2 \indi{  |V_n| > \frac 1 {\eta \sqrt{n \log n}} }\right] = 0.
\end{equation}
Let $n$ be so large that $\frac{1}{\eta \sqrt{n \log n}} < R$. We split the above  expectations up into
\begin{align}
\E \left[Z_n \overline{Z_n} \indi{  |V_n| > \frac 1 {\eta \sqrt{n \log n}} }\right]
&=
\E \left[Z_n \overline{Z_n} \indi{ \frac 1 {\eta \sqrt{n \log n}} < |V_n| < R }\right]
+\E \left[Z_n \overline{Z_n} \indi{ |V_n| \ge R  }\right], \label{eq:exp_split_up1}\\
\E \left[Z_n^2 \indi{ |V_n| > \frac 1 {\eta \sqrt{n \log n}} }\right]&=
\E \left[ Z_n^2 \indi{ \frac 1 {\eta \sqrt{n \log n}} < |V_n| < R }\right]
+\E \left[Z_n^2 \indi{ |V_n| \ge R  }\right]. \label{eq:exp_split_up2}
\end{align}
For the second summands on the right-hand sides of~\eqref{eq:exp_split_up1} and~\eqref{eq:exp_split_up2} we have the estimates
\begin{align*}
\E \left[Z_n \overline{Z_n} \indi{ |V_n| \ge R }\right]
&=\frac{1}{n \log n} \E \left[\frac{1}{|V_n|^2} \indi{ |V_n| \ge R  }\right]
\le \frac{1}{R^2 n \log n} = o \left( \frac{1}{n} \right),\\
\left| \E  Z_n^2 \indi{ V_n \ge R } \right|
&\leq
\E \left[Z_n \overline{Z_n} \indi{ |V_n| \ge R }\right]
=o \left(\frac{1}{n} \right),
\end{align*}
which gives
\begin{align}
\lim_{n\to\infty} n \E \left[Z_n \overline{Z_n} \indi{ V_n \ge R }\right] = 0
\quad
\text{and}
\quad
\lim_{n\to\infty}  n \E \left[ Z_n^2 \indi{ V_n \ge R } \right] = 0.
\end{align}
We now analyze the first summands on the right-hand sides of~\eqref{eq:exp_split_up1} and~\eqref{eq:exp_split_up2}. To this end, we shall use Lemma~\ref{lemma:Integralasymptote}.
Recalling that the density of $V_n$ near $0$ is $q_n$ and passing to  polar coordinates, we can write
\begin{align*}
n \E \left[Z_n \overline{Z_n} \indi{ \frac{1}{\eta \sqrt{n \log n}} < \left| V_n \right| < R }\right]
&=
\frac{1}{\log n} \int_{ \left\{ \frac{1}{\eta \sqrt{n \log n}} < \left| v \right| < R \right\} } \frac{q_n(v)}{|v|^2}  \dd v\\
&=
\frac{1}{\log n} \int_\frac{1}{\eta \sqrt{n \log n}}^{R} \frac{1}{r} \int_0^{2 \pi} q_n(r\eee^{i\phi}) \dd \phi \dd r.
\end{align*}
Here we can use Lemma \ref{lemma:Integralasymptote} with $\delta(n) : =\frac{1}{\eta \sqrt{n \log n}}$,
$g_{\delta(n)} (r) = \int_0^{2 \pi} q_{n}(r\eee^{i\phi}) \dd \phi$ and $g_{0} (r) = \int_0^{2 \pi} p(r\eee^{i\phi}) \dd \phi$. Note that $g_0(0)= 2\pi p(0)$.
The assumptions of Lemma~\ref{lemma:Integralasymptote} are fulfilled since by the continuity of $p$ at $0$ and the condition $\lim_{n\to\infty} z_n = 0$, for every $\eps>0$ we have
$$
|g_{\delta(n)} (r)-g_0(0)| = \left|\int_0^{2 \pi} (q_{n}(r\eee^{i\phi})-p(0)) \dd \phi\right|
\leq \int_0^{2 \pi} |p(z_n-r\eee^{i\phi})-p(0)| \dd \phi \leq \eps
$$
if $n$ is sufficiently large and $r$ is sufficiently small.
Applying Lemma~\ref{lemma:Integralasymptote} we obtain
$$
n \left[\E Z_n \overline{Z_n} \indi{ \frac{1}{\eta \sqrt{n \log n}} < \left| V_n \right| < R }\right]
=
\frac{\log \left( \eta \sqrt{n \log n} \right)}{\log n}  \frac{1}{\log \frac{1}{\delta(n)}} \int_{\delta(n)}^{R} \frac{g_{\delta(n)} (r)}{r} \dd r
\ton
\frac{g_0(0)}{2}
=\pi p(0).
$$
Similarly, using Lemma \ref{lemma:Integralasymptote} with $\delta(n) : =\frac{1}{\eta \sqrt{n \log n}}$, $g_{\delta(n)} (r) = \int_0^{2 \pi} \eee^{-2i \phi} q_{n}(r\eee^{i\phi}) \dd \phi$ and $g_{0} (r) = \int_0^{2 \pi} \eee^{-2i \phi} p(r\eee^{i\phi}) \dd \phi$, we obtain
\begin{multline*}
n \E \left[Z_n^2 \indi{ \frac{1}{\eta \sqrt{n \log n}} < \left| V_n \right| < R }\right]
=
\frac{1}{\log n} \int_{\left\{ \frac{1}{\eta \sqrt{n \log n}} < \left| v \right| < R \right\}} \frac{q_n(v)}{v^2} \dd v  \\
=
\frac{1}{\log n} \int_{\frac{1}{\eta \sqrt{n \log n}}}^{R} \frac{1}{r} \int_0^{2 \pi} \eee^{-2i \phi} q_n(r\eee^{i \phi}) \dd \phi \dd r
=
\frac{\log \left( \eta \sqrt{n \log n} \right)}{\log n} \frac{1}{\log \frac 1{\delta(n)}} \int_{\delta(n)}^{R} \frac{g_{\delta(n)}(r)}{r}  \dd r
\ton 0
\end{multline*}
because $g_0(0) = \int_0^{2 \pi} \eee^{-2i \phi} p(0) \dd \phi =0$. This completes the proof of~\eqref{eq:exp_z_n_overline_z_n} and thus of~\eqref{eq:cov_lim1} and~\eqref{eq:cov_lim2}.

\medskip
\noindent
\textit{Proof of~\eqref{eq:cov_lim3}.}
Since both $|\E [\Re Z_n \indi{ \left| Z_n \right| < \eta }]|$ and $|\E [\Im Z_n \indi{ \left| Z_n \right| < \eta }]|$ can be upper bounded by  $\E | Z_n|$, it suffices to prove that
\begin{equation}\label{eq:lim_sqrt_n_E_Z_n_0}
\lim_{n\to\infty} \sqrt n \, \E |Z_n| = 0.
\end{equation}
Recalling that $Z_n = 1/(V_n \sqrt{n\log n})$ and that on the disk $\bB_R(0)$, the random variable $V_n$ has density $q_n$ bounded by $C$, we can write
\begin{align*}
\E |Z_n|
&=
\frac{1}{\sqrt{n \log n}} \E \frac 1 {|V_n|}
=
\frac{1}{\sqrt{n \log n}} \int_{\bB_R(0)} \frac{q_n(v)}{|v|} \dd v
+
\frac{1}{\sqrt{n \log n}} \E \left[\frac 1 {|V_n|} \ind_{\{|V_n| > R\}} \right]\\
&\leq
\frac{1}{\sqrt{n \log n}} \int_{\bB_R(0)} \frac{C}{|v|} \dd v
+
\frac{1}{\sqrt{n \log n}} \frac 1 {R}
=
O\left(\frac{1}{\sqrt{n \log n}}\right),
\end{align*}
which proves~\eqref{eq:lim_sqrt_n_E_Z_n_0}.
The proof of~\eqref{eq:meerschaerd2_applied} is thus complete.

\medskip
Now we can apply Theorem \ref{theo:rvaceva} that yields, for every $\rho>0$,
\begin{align}
\frac{1}{\sqrt{n \log n}} \sum_{k=1}^n \left( \frac{1}{z_n-\xi_k} - \E \left[ \frac{1}{z_n-\xi_k} \indi{ \left| \frac{1}{z_n-\xi_k} \right| < \rho \sqrt{n \log n}} \right] \right) \todistr  \cN_\C \left( 0 , \pi p(0) \right).
\end{align}
In view of Slutsky's lemma, to complete the proof of Theorem~\ref{theo:CLT}, we need to show that
\begin{align}
\lim_{n\to\infty} \sqrt{ \frac{n}{\log n} } \left( \E \left[ \frac{1}{z_n-\xi_k} \indi{ \left| \frac{1}{z_n-\xi_k} \right| < \rho \sqrt{n \log n} }\right] -f(0) \right) = 0. \label{eq:VorSlutsky}
\end{align}

\medskip
\noindent
\textit{Proof of~\eqref{eq:VorSlutsky}.}
Since the density $p$ of the $\xi_k$'s  exists and is bounded in a neighborhood of $0$, and since $\lim_{n\to\infty} z_n = 0$, for sufficiently large $n$ we have the estimate
\begin{align}
\left| \E \left[ \frac{1}{z_n-\xi_k} \indi{\left| \frac{1}{z_n-\xi_k} \right| \ge \rho \sqrt{n \log n}} \right] \right|
&\leq
\int_{\bB_{\frac{1}{\rho \sqrt{n \log n}}}(z_n)}  \frac{p(u)}{|z_n-u|}  \dd u
\leq
C \int_{\bB_{\frac{1}{\rho \sqrt{n \log n}}}(z_n)} \left| \frac{1}{u-z_n} \right| \dd u\notag\\
&=
C \int_{\bB_{\frac{1}{\rho \sqrt{n \log n}}}(0)} \frac{1}{|v|} \dd v
=
\frac{2 \pi C}{\sqrt{n \log n} \rho}
=
O\left( \frac{1}{\sqrt{n \log n}} \right). \label{eq:proof_of_34}
\end{align}
Thus,
\begin{align*}
\E \left[ \frac{1}{z_n-\xi_k} \indi{ \left| \frac{1}{z_n-\xi_k} \right| < \rho \sqrt{n \log n} } \right]
= \E \left[ \frac{1}{z_n-\xi_k}\right] + O\left( \frac{1}{\sqrt{n \log n}} \right)
=f(z_n) +O \left( \frac{1}{\sqrt{n \log n}} \right).
\end{align*}
Using Lemma~\ref{lem:lipschitz}, below, together with the monotone increasing property of the function $x \mapsto |x \log x|$, $0<x<1/\eee$, and the condition $z_n=O(n^{-1/2-\kappa})$ with $\kappa>0$, we obtain
\begin{align*}
\E \left[ \frac{1}{z_n-\xi_k} \indi{\left| \frac{1}{z_n-\xi_k} \right| < \rho \sqrt{n \log n} } \right] -f(0)
&= f(z_n) -f(0) +O \left( \frac{1}{\sqrt{n \log n}} \right)\\
&= O(|z_n \log |z_n||) +O \left( \frac{1}{\sqrt{n \log n}} \right)
=O \left( \frac{1}{\sqrt{n \log n}} \right) \text{,}
\end{align*}
thus proving~\eqref{eq:VorSlutsky}.
 This completes the proof of Theorem~\ref{theo:CLT}.
\end{proof}

\begin{remark}\label{rem:after_34}
For later use, observe that~\eqref{eq:VorSlutsky} continues to hold if $\rho \sqrt{n\log n}$ is replaced by any larger sequence $d_n$. Indeed, after this replacement~\eqref{eq:proof_of_34} holds with the better error estimate $O(1/d_n)$, while the rest of the proof remains the same.
\end{remark}

The following lemma, which we already used above, will be essential at several places in the proof.
\begin{lemma}\label{lem:lipschitz}
Let $U$ be a random variable that has a Lebesgue density $p$ on some disk $\bB_r(0)$ and may have arbitrary distribution outside this disk. Assume also that $p$ is bounded by a constant $c_1$ on $\bB_r(0)$. Then, the Cauchy--Stieltjes transform $f(z):=\E [\frac{1}{z-U}]$ exists finitely on $\bB_{r/2}(0)$ and for a suitable constant $C>0$ we have
\begin{equation}\label{eq:lem:lipschitz}
|f(z) - f(0)| \leq C |z \log |z||
\quad \text{ for all } z\in\C \text{ such that } |z|< \frac 12 \min\{1,r\}.
\end{equation}
\end{lemma}

\begin{proof}
To prove the finiteness of the Cauchy--Stieltjes transform on $\bB_{r/2}(0)$, take some $z\in \bB_{r/2}(0)$ and write
$$
\E \left|\frac 1 {z-U} \right| = \E \left[\left|\frac 1 {z-U} \right|\ind_{\{U\in \bB_{r/2}(z)\}} \right] + \E \left[\left|\frac 1 {z-U} \right|\ind_{\{U\notin \bB_{r/2}(z)\}} \right].
$$
The second expectation can be trivially bounded by $2/r$, so let us consider the first one. The density of the random variable $z-U$, denoted by $q(x)$, exists on the ball $\bB_{r/2}(0)$ and is bounded by $c_1$ there. Hence,
$$
\E \left[\left|\frac 1 {z-U} \right|\ind_{\{U\in \bB_{r/2}(z)\}} \right]
= \int_{\bB_{r/2}(0)} \frac{q(w)}{|w|}  \dd w
\leq
c_1 \int_{\bB_{r/2}(0)} \frac{\dd w}{|w|}
=
c_1 \int_{0}^{r/2} \int_0^{2\pi}  \dd s\dd\theta
=
\pi r c_1,
$$
which is finite.
Let us prove~\eqref{eq:lem:lipschitz}.  By definition of $f$, we have
$$
|f(z) - f(0)| = \left|\E \left[ \frac 1 {z-U} + \frac 1 U\right]\right|
=
\left| \E \left[\frac {z}{(z-U)U}\right]\right|
\leq
|z|   \E \left[\frac {1}{|z-U||U|}\right].
$$
In the following, let $z\in\C$ be such that $|z|< \frac 12 \min\{1,r\}$. Our aim is to show that
\begin{equation}\label{eq:f_lipschitz_aim}
\E \left[\frac {1}{|z-U||U|}\right] \leq C |\log |z||.
\end{equation}
First of all, if $|U|\geq r$, then by the triangle inequality, $|z-U|\geq \frac r2$ and hence
$$
\E \left[\frac {1}{|z-U||U|} \ind_{\{|U|\geq r\}}\right] \leq  \E \left[\frac 1{r^2/2}\right] = \frac 2{r^2} = \text{const}
< \frac C4 |\log |z||
$$
provided $C$ is sufficiently large.
If $2 |z| \leq |U| < r$, then $|U-z|\geq |U|-|z| \geq \frac 12 |U|$ by the triangle inequality and therefore
\begin{multline*}
\E \left[\frac {1}{|z-U||U|} \ind_{\{2|z| \leq |U|<  r\}}\right]
\leq
\E \left[\frac {2}{|U|^2} \ind_{\{2|z| \leq |U|<  r\}}\right]
=
\int_{2|z|}^r \int_0^{2\pi} \frac{2}{s^2} p (s\eee^{i \theta}) s \dd s\dd\theta
\\
\leq
4\pi \int_{2|z|}^r \frac {c_1} {s} \dd s
=
4\pi c_1 (\log r - \log (2|z|))
=
\text{const} - 4\pi c_1 |\log |z| |
\leq
\frac C4 |\log |z||.
\end{multline*}
It remains to estimate the expectation on the event $F= \{|U| \leq  2 |z|\}$. Consider the events $E_1 = \{|U| \leq |U-z|\}$ and $E_2 = \{|U| \geq |U-z|\}$. By the triangle inequality,  $|U-z| + |U| \geq  |z|$, hence on the event $F\cap E_1$ we have $|U-z| \geq \frac 12 |z|$ and
\begin{multline*}
\E \left[\frac {1}{|z-U||U|} \ind_{F\cap E_1}\right]
\leq
\E \left[\frac {2}{|U||z|} \ind_{F}\right]
=
\int_{0}^{2|z|} \int_0^{2\pi} \frac{2}{r|z|} p (r\eee^{i \theta}) r \dd r\dd\theta
\\
\leq
\frac {4\pi}{|z|} \int_{0}^{2|z|}c_1 \dd r
\leq 8\pi c_1
=\text{const}
\leq
\frac C4 |\log |z||.
\end{multline*}
Similarly, on the event $F\cap E_2$ we have $|U| \geq \frac 12 |z|$. Thus,
\begin{multline*}
\E \left[\frac {1}{|z-U||U|} \ind_{F\cap E_2}\right]
\leq
\E \left[\frac {2}{|U-z||z|} \ind_{F}\right]
=
\int_{\bB_{2|z|}(0)} \frac{2 p (w) \dd w}{|w-z||z|}
\leq
\frac {2c_1}{|z|} \int_{\bB_{2|z|}(0)} \frac{\dd w}{|w-z|}
\\
\leq
\frac {2c_1}{|z|} \int_{\bB_{3|z|}(z)} \frac{\dd w}{|w-z|}
=
\frac {2c_1}{|z|}  \int_{0}^{3|z|} \int_0^{2\pi} \frac{1}{r} r \dd r\dd\theta
=
\text{const}
\leq
\frac C4 |\log |z||,
\end{multline*}
where we passed to polar coordinates with the origin shifted to the point $z$.
Taking everything together, we arrive at~\eqref{eq:f_lipschitz_aim}, thus completing the proof.
\end{proof}

\subsection{Proof of Theorem~\ref{theo:funct.CLT}}
We divide the proof into  two parts. First we show that the finite-dimensional distributions converge and then we shall prove tightness.

\begin{lemma} \label{lemma:funct.CLT-fdd}
Under the assumptions of Theorem~\ref{theo:funct.CLT}, for all $d\in\N$ and $w_1,\ldots,w_d\in \C$, the following weak convergence of random vectors holds true:
\begin{align*}
\left(\frac{1}{\sqrt{n \log n}} \sum_{k=1}^n \left( \frac{1}{z_n(w_j)-\xi_k} - f(0) \right)\right)_{j=1,\ldots,d} \todistr
(N,\ldots,N),
\end{align*}
where $N\sim \cN_\C ( 0 , \pi p(0))$.  The components of the limit vector almost surely are equal.
\end{lemma}
\begin{proof}
By Theorem~\ref{theo:CLT}, we have
$$
\frac{1}{\sqrt{n \log n}} \sum_{k=1}^n \left( \frac{1}{z_n(w_1)-\xi_k} - f(0) \right) \todistr
N.
$$
To prove the lemma it is sufficient to show that for all $i\in\{2,\ldots,d\}$,
\begin{align}
\left|
\frac{1}{\sqrt{n \log n}} \sum_{k=1}^n \left( \frac{1}{z_n(w_1)-\xi_k} - f(0) \right)
-
\frac{1}{\sqrt{n \log n}} \sum_{k=1}^n \left( \frac{1}{z_n(w_i)-\xi_k} - f(0) \right)
\right|
\toprobab 0 .
\end{align}
Using the definition of $z_n(w)$ we obtain
\begin{align*}
&
\left| \frac{1}{\sqrt{n \log n}} \sum_{k=1}^n \left( \frac{1}{z_n(w_1)-\xi_k} - f(0) \right)
-
\frac{1}{\sqrt{n \log n}} \sum_{k=1}^n \left( \frac{1}{z_n(w_i)-\xi_k} - f(0) \right) \right| \\
&=
\frac{1}{\sqrt{n \log n}} \left| \sum_{k=1}^n \frac{z_n(w_i)-z_n(w_1)}{(z_n(w_1)-\xi_k)(z_n(w_i)-\xi_k)} \right|\\
&=
\frac{1}{n^2} \left| \sum_{k=1}^n \frac{ w_i-w_1}{(z_n(w_1)-\xi_k)(z_n(w_i)-\xi_k)} \right| \cdot \frac 1 {|f(0)|^2}\\
&\leq
\frac{|w_i-w_1|}{n^2} \left| \sup_{\zeta_1,\zeta_2\in\C\colon |\zeta_1|, |\zeta_2| \le \frac{C}{n}}
\sum_{k=1}^n \frac{1}{(\zeta_1-\xi_k)(\zeta_2-\xi_k)} \right| \cdot \frac 1 {|f(0)|^2}
\end{align*}
since $|z_n(w_1)|$ and $|z_n(w_i)|$ can be upper bounded by $\frac Cn$ for some constant $C$.
The proof is completed by a use of Lemma \ref{lemma:convergence sum of squares}, see below, with $s_n=\frac{C}{n}$ and $a_n =\frac{1}{n^2}$.
\end{proof}

\begin{lemma} \label{lemma:convergence sum of squares}
In the situation of Theorem~\ref{theo:funct.CLT},  fix  positive sequences $s_n \in (0,1)$ and $a_n>0$ satisfying
\begin{equation}\label{eq:a_n_s_n_requirements}
s_n = o \left( \frac{1}{\sqrt{n}} \right), \quad a_n = o \left(\frac{1}{n |\log s_n|} \right),  \quad n\to\infty.
\end{equation}
Writing $\sup_{|z_1|, |z_2| \le s_n}$ as a shorthand for  $\sup_{z_1, z_2 \in \C :|z_1|, |z_2| \le s_n}$, we have
\begin{align}
a_n \sup_{|z_1|, |z_2| \le s_n} \left| \sum_{k=1}^n \frac{1}{\left(z_1-\xi_k \right) \left(z_2-\xi_k \right)} \right| \toprobab 0 \text{.}
\end{align}
\end{lemma}
\begin{proof}
Consider the random event $A_n := \left\{|\xi_1| > 2 s_n, \ldots, |\xi_n| > 2 s_n \right\}$. It suffices to show that
\begin{equation}\label{eq:exp_sup_z1_z2_on_An}
\lim_{n\to\infty} a_n \E \left[ \sup_{|z_1|, |z_2| \le s_n} \left| \sum_{k=1}^n \frac{1}{\left(z_1-\xi_k \right) \left(z_2-\xi_k \right)} \right|\ind_{A_n}\right] = 0
\end{equation}
and $\lim_{n\to\infty} \P[A_n^c] = 0$.
By the triangle inequality, on the event $A_n$ we have $|z_j-\xi_k| \ge \frac{|\xi_k|}{2}$ for all $j \in \{ 1,2 \}$ and $k\in\{1,\ldots,n\}$. Using this together with the trivial estimate $\ind_{A_n} \leq \indi{|\xi_k| >2s_n }$, we arrive at
\begin{align*}
a_n \E \left[\sup_{|z_1|, |z_2| \le s_n} \left| \sum_{k=1}^n \frac{1}{\left(z_1-\xi_k \right) \left(z_2-\xi_k \right)}  \right|\ind_{A_n}\right]
&\le
a_n  \E \left[\sum_{k=1}^n   \frac{4}{|\xi_k|^2} \indi{|\xi_k|>2s_n}\right]\\
&=
4 a_n n  \E \left[\frac{1}{|\xi_1|^2} \indi{|\xi_1|>2s_n}\right].
\end{align*}
For a sufficiently small $r>0$, the density $p$ of $\xi_1$ satisfies $\sup_{z\in \bB_r(0)} |p(z)| \leq C$ by the continuity of $p$ at $0$.
Splitting the expectation and passing to the polar coordinate system, we obtain
\begin{multline*}
4 a_n n  \E \left[\frac{1}{|\xi_1|^2} \indi{|\xi_1|>2s_n}\right]
=
4 a_n n \int_{\bB_r(0) \setminus \bB_{2s_n}(0)} \frac{p(u)}{|u|^2} \dd u  +
4 a_n n \E \left[\frac{1}{|\xi_1|^2} \indi{|\xi_1|>r}\right]\\
\le 4 a_n n \left(2 \pi C\int_{2s_n}^r \frac{1}{s} \dd s +O(1) \right)
= 4 a_n n \left( -2 \pi C  \log (2 s_n) + O(1) \right) \text{,}
\end{multline*}
Under the given requirements~\eqref{eq:a_n_s_n_requirements}, this upper bound goes to $0$, thus proving~\eqref{eq:exp_sup_z1_z2_on_An}.

It remains to prove that $\lim_{n\to\infty} \P[A_n^c] = 0$.
Since the density $p$ exists and is bounded on a small neighborhood of $0$ and since $\lim_{n\to\infty} s_n = 0$, for sufficiently large $n$ the union bound yields
\begin{align}
\P[A_n^c] = \P[\exists k\in \{1,\ldots,n\}:|\xi_k| \leq 2s_n ] \leq n \P[|\xi_1| \leq 2s_n] \leq n 4\pi s_n^2 C, \label{eq:convergence of squares estimate of probability}
\end{align}
which converges to $0$ since $s_n^2= o( \frac{1}{n})$  by~\eqref{eq:a_n_s_n_requirements}.
\end{proof}

In the proof of Theorem~\ref{theo:funct.CLT} the remaining step is to prove tightness.
\begin{lemma} \label{lemma:funct.CLT-tightness}
Under the assumptions of Theorem~\ref{theo:funct.CLT}, the sequence of stochastic processes
\begin{align}
\left(\frac{1}{\sqrt{n \log n}} \sum_{k=1}^n \left( \frac{1}{z_n(w) - \xi_k} - f(0) \right)\right)_{|w|\leq \rho}, \qquad n\geq n_0,
\end{align}
is tight on the Banach space $\bA_{\rho}$.
\end{lemma}
\begin{proof}
We split the sequence into
\begin{align}
\frac{1}{\sqrt{n \log n}} \sum_{k=1}^n &\left( \frac{1}{z_n(w)-\xi_k} - f(0) \right) \\
&=\frac{1}{\sqrt{n \log n}} \sum_{k=1}^n \left( \frac{1}{z_n(w)-\xi_k} -f(0) \right) \left( \indi{|\xi_k| < n^{-\frac{2}{3}}} + \indi{|\xi_k| \ge n^{-\frac{2}{3}}} \right).
\end{align}
In fact,  instead of $n^{-\frac{2}{3}}$ we could have chosen any sequence $a_n>0$ satisfying
$$
\lim_{n \to \infty} na_n =+\infty, \quad
\lim_{n \to \infty} \sqrt n a_n =0.
$$
It has already been proven in~\eqref{eq:convergence of squares estimate of probability}, in the proof on Lemma \ref{lemma:convergence sum of squares}, that
\begin{align*}
\P \left[ \exists k \in \{ 1,\ldots,n \} : \left| \xi_k \right| < n^{- \frac{2}{3}} \right] \ton 0.
\end{align*}
This  implies that
\begin{align}
\bP \left[ \frac{1}{\sqrt{n \log n}} \sum_{k=1}^n \left( \frac{1}{z_n(\cdot)-\xi_k} -f(0) \right) \indi{|\xi_k| < n^{-\frac{2}{3}}} \equiv 0 \right] \ton 1. \label{eq:proof_funct.CLT_cond2}
\end{align}
Therefore,  it suffices to show that the sequence of stochastic processes
\begin{align}
\left(\frac{1}{\sqrt{n \log n}} \sum_{k=1}^n \left( \frac{1}{z_n(w)-\xi_k} -f(0) \right) \indi{|\xi_k| \ge n^{-\frac{2}{3}}}\right)_{|w|\leq \rho},
\qquad n\geq n_0,
\label{eq:proof_funct.CLT_cond1}
\end{align}
is tight on $\bA_\rho$. 
%
To this end, we shall show that
\begin{align*}
\beta_n(w) := \E \left| \frac{1}{\sqrt{n \log n}} \sum_{k=1}^n \left( \frac{1}{z_n(w)-\xi_k} -f(0) \right) \indi{|\xi_k| \ge n^{-\frac{2}{3}}} \right|^2
\end{align*}
is bounded by a constant that does not depend on $n\geq n_0$ and $|w|\leq \rho$. This implies tightness on $\bA_{\rho}$ by~\cite[Remark on p. 341]{shirai}.
Since the summands are i.i.d.\ (but their mean is, in general, non-zero), we have
\begin{multline*}
\beta_n(w)
=
\frac{1}{\log n} \E \left[ \left| \frac{1}{z_n(w)-\xi_1} -f(0) \right|^2 \indi{|\xi_1| \ge n^{-\frac{2}{3}}} \right]
\\+ \frac{n-1}{\log n} \left|\E \left[ \left(\frac{1}{z_n(w)-\xi_1} -f(0)\right)\indi{|\xi_1| \ge n^{-\frac{2}{3}}} \right] \right|^2
\end{multline*}
The second summand goes to $0$ by~\eqref{eq:VorSlutsky}; see also Remark~\ref{rem:after_34}.
In fact, this convergence is even uniform in $w$ as long as $|w| \leq \rho$ since~\eqref{eq:VorSlutsky} holds uniformly in $|z_n| = O(n^{-1/2-\kappa})$, which is fulfilled for $z_n= z_n(w)$.
Turning to  the first summand, the inequality $|a+b|^2 \leq 2|a|^2 + 2|b|^2$ leads to the estimate
$$
\frac{1}{\log n} \E \left[ \left| \frac{1}{z_n(w)-\xi_1} -f(0) \right|^2 \indi{|\xi_1| \ge n^{-\frac{2}{3}}} \right]
\leq
\frac{2}{\log n} \E \left[ \left| \frac{1}{z_n(w)-\xi_1}\right|^2 \indi{|\xi_1| \ge n^{-\frac{2}{3}}} \right] + \frac{2|f(0)|^2}{\log n}.
$$
The second term is $O(1)$, so let us consider the first one. The subsequent estimates are uniform over $|w|\leq \rho$.
Splitting the expectation and using the inequality $|\xi_1-z_n(w)| \ge \frac{|\xi_1|}{2}$, which holds for large $n$ on the event $|\xi_1| \ge n^{-\frac{2}{3}}$ because $z_n(w) = O(\frac 1n)$, we obtain
\begin{multline*}
\frac{2}{\log n} \E \left[ \left| \frac{1}{z_n(w)-\xi_1}\right|^2 \indi{|\xi_1| \ge n^{-\frac{2}{3}}} \right]
\\
\leq
\frac{2}{\log n} \E \left[ \left| \frac{1}{z_n(w)-\xi_1}\right|^2 \indi{|\xi_1| \ge R} \right]
+
\frac{8}{\log n} \E \left[ \frac{1}{|\xi_1|^2} \indi{ n^{-\frac{2}{3}} \leq |\xi_1| < R} \right]
\end{multline*}
The first expectation on the right-hand side is $O(1)$ because $|z_n(w)-\xi_1|$ is bounded away from $0$ on the event $\{|\xi_1| \geq R\}$. To bound the second summand,  we pass to polar coordinates and use the fact that the density $p$ exists and is bounded by a constant $C$ on the disk $\bB_{R}(0)$, thus arriving at
\begin{multline*}
\frac{8}{\log n} \E \left[ \frac{1}{|\xi_1|^2} \indi{ n^{-\frac{2}{3}} \leq |\xi_1| < R} \right]
=
\frac{8}{\log n} \int_{n^{-\frac{2}{3}}}^R \int_0^{2 \pi} \frac{p(r \eee^{i\phi})}{r} \dd \phi \dd r
\\
\leq
\frac{16 \pi C}{\log n} \int_{n^{-\frac{2}{3}}}^{R} \frac{1}{r} \dd r
=
\frac{32}{3} \pi C + O\left(\frac 1{\log n}\right)
=
O(1).
\end{multline*}
Taking everything together, we obtain $\beta_n(w) = O(1)$ uniformly over $n\geq n_0$ and $|w|\leq \rho$,  thus proving that the sequence  in~\eqref{eq:proof_funct.CLT_cond1} is tight.
\end{proof}

\subsection{Proof of Theorem~\ref{theo:criticalpoints}}
Without restriction of generality we assume that $u_0=0$.
We shall study the zeros of the logarithmic derivative
\begin{equation}\label{eq:P_n_log_der}
\frac{P_n'(z)}{P_n(z)}
=
\frac{\dd}{\dd z} \log P_n(z)
=
\frac{\dd}{\dd z} \left( \log z +\sum_{k=1}^n \log (z-\xi_k) \right)
=
\frac{1}{z} +\sum_{k=1}^n \frac{1}{z-\xi_k}
\end{equation}
instead of those of $P_n'$. Indeed, since the Lebesgue density $p$ of the random variables $\xi_k$ exists in some neighborhood of $0$, the probability that $P_n$, for some $n$,  has a multiple zero in this neighborhood is $0$. Thus, for sufficiently large $n$, the polynomial $P_n$ has no multiple zeroes in $\bB_{r_n}(0)$, with probability $1$, and hence the zeroes of $P_n'/P_n$ in $\bB_{r_n}(0)$ coincide with those of $P_n$.

\medskip
\noindent
\textit{Proof of~\eqref{eq:mainthm1}.}
Our aim is to show that with probability converging to $1$, the polynomial $P_n$ has exactly one critical point in the disk $\bB_{r_n}(0)$.
Let $\bT_r(0)=\{z\in\C\colon |z| = r\}$ be the circle with radius $r>0$ that is centered at the origin. 
The main step in proving \eqref{eq:mainthm1} is the following
\begin{lemma}\label{lem:proofRouche first equation}
We have
\begin{align}
\sup_{z\in \bT_{r_n}(0)}
\left|
\frac{1}{n} \frac{P_n'(z)}{P_n(z)} - f(0) \right|
\toprobab 0. \label{eq:proofRouche first equation}
\end{align}
\end{lemma}
\begin{proof}
First of all, observe that by~\eqref{eq:P_n_log_der},
$$
\sup_{z\in \bT_{r_n}(0)}
\left|
\frac{1}{n} \frac{P_n'(z)}{P_n(z)} - f(0) \right|
=
\sup_{z\in \bT_{r_n}(0)}
\left|
\frac{1}{nz} + \frac 1n \sum_{k=1}^n \frac{1}{z-\xi_k} - f(0)
\right|.
$$
Since
$$
\lim_{n\to\infty} \sup_{z\in \bT_{r_n}(0)} \left|\frac 1 {nz}\right| = \lim_{n\to\infty} \frac 1 {nr_n} = 0,
$$
to prove the lemma it suffices to show that
\begin{equation}\label{eq:proofRouche first equation_repeat}
\sup_{z \in \bB_{r_n}(0)} \left|\frac{1}{n} \sum_{k=1}^n \frac{1}{z_n-\xi_k} - f(0)\right| \toprobab 0.
\end{equation}
As a first step we shall prove the following weak law of large numbers: For every complex sequence $(z_n)_{n\in\N}$ converging to $0$, we have
\begin{align}
\frac{1}{n} \sum_{k=1}^n \frac{1}{z_n-\xi_k} \toprobab f(0). \label{eq:WLLN_second}
\end{align}
Fix  some $\eps >0$ and $p \in (1,2)$.
An inequality due to von Bahr and Esseen~\cite{von_bahr_esseen} states that for
centered, independent real or complex random
variables $\eta_1,\ldots,\eta_n$ with finite $p$-th absolute moment we have
\begin{align} \label{eq:ineq_moment_p_02}
\E |\eta_1+\ldots+\eta_n|^p \leq 2 \sum_{k=1}^n \E |\eta_k|^p.
\end{align}
Using the inequalities of Markov, von Bahr--Esseen, and finally the inequality $|a+b|^p \leq 2^{p-1}(|a|^p + |b|^p)$, we obtain
\begin{multline}\label{eq:von_bahr}
\P \left[ \left| \frac{1}{n} \sum_{k=1}^n \frac{1}{z_n-\xi_k} -f(z_n) \right| > \eps \right]
\leq
\frac{1} {(n \eps)^p}
\E \left| \sum_{k=1}^n \left( \frac{1}{z_n-\xi_k} -f(z_n) \right) \right|^p
\\
\leq
\frac{2}{(n \eps)^p} \sum _{k=1}^n \E \left| \frac{1}{z_n-\xi_k}-f(z_n) \right|^p
\leq
\frac{2^{p}}{n^{p-1} \eps^p}  \left( \E \left| \frac{1}{z_n-\xi_1} \right|^p +\left|f(z_n) \right|^p \right).
\end{multline}
We have $\left|f(z_n) \right|^p=O(1)$ since in fact $\lim_{n\to\infty} f(z_n) = f(0)$ by Lemma~\ref{lem:lipschitz}. Let $q_n(z) =p(z_n - z)$ be the density of $z_n-\xi_1$.
Since $p$ is continuous at $0$, there is a sufficiently small $R>0$ such that $p$ is bounded in the disk $\bB_{2R}(0)$ by some constant $C$.
It follows from $\lim_{n\to \infty} z_n = 0$ that  the density  $q_n$ is bounded on the smaller disk $\bB_{R}(0)$ by the same constant $C$  provided $n$ is sufficiently large. Thus,
$$
\E \left|\frac{1}{z_n-\xi_1} \right|^p
=
\E \left[\left|\frac{1}{z_n-\xi_1} \right|^p \ind_{\{|z_n-\xi_1| > R\}}\right]
+
\int_{\bB_R(0)} \frac{q_n(z)}{|z|^p}\dd z
\leq
\frac1 {R^p} +
\int_{\bB_R(0)} \frac{C}{|z|^p}\dd z
=
O(1).
$$
Hence, the probability on the right-hand side of~\eqref{eq:von_bahr} converges to $0$ proving that
\begin{align}
\frac{1}{n} \sum_{k=1}^n \frac{1}{z_n-\xi_k}-f(z_n) \toprobab 0. \label{eq:WLLN}
\end{align}
To complete the proof of~\eqref{eq:WLLN_second} recall that $\lim_{n\to\infty} f(z_n) = f(0)$ by Lemma~\ref{lem:lipschitz}.
We prove~\eqref{eq:proofRouche first equation_repeat}.
By Lemma~\ref{lemma:convergence sum of squares} with $s_n = r_n$ and $a_n = \frac 1n$ we have
\begin{align}
\sup_{z', z''\in \bB_{r_n}(0)} \left| \frac{1}{n} \sum_{k=1}^n \left( \frac{1}{z'-\xi_k} -\frac{1}{z''-\xi_k} \right) \right| \toprobab 0 \text{.} \label{eq:WLLN uniformly}
\end{align}
Together with~\eqref{eq:WLLN_second} this yields~\eqref{eq:proofRouche first equation_repeat} by the triangle inequality, thus completing the proof of the lemma.
\end{proof}


Now we are in position to prove~\eqref{eq:mainthm1}.
By Rouch\'e's theorem~\cite[pp.125--126]{conway_book}, if the event
$$
E_n^{(1)} :=  \left\{\sup_{z\in \bT_{r_n}(0)}\left| \frac{1}{n} \frac{P_n'(z)}{P_n(z)}-f(0)\right| < |f(0)|\right\}
$$
occurs, then the difference between the number of zeroes and the number of poles of the function $\frac{1}{n} \frac{P_n'}{P_n}$ on the disk $\bB_{r_n}(0)$ is equal to the difference between the numbers of zeroes and poles of the constant function $z \mapsto f(0)$, which is $0$. Since $f(0) \neq 0$, \eqref{eq:proofRouche first equation} implies that
\begin{align}
\lim_{n\to\infty}
\P \left[ E_n^{(1)}\right] = 1.
\label{eq:convergence for Rouche}
\end{align}
On the other hand, the probability of the event
$$
E_n^{(2)}:= \left\{ \frac{P_n'}{P_n} \text{ has exactly one pole in } \bB_{r_n}(0)\right\}
$$
converges to $1$ since $0$ is a pole of $P_n'/P_n$ and thus
\begin{align*}
\P \left[E_n^{(2)}\right]
=
\P [\text{there is no } z \in \bB_{r_n}(0) \backslash\{0\} \colon P_n(z)=0]
=
\P \left[ \bigcap_{k=1}^n \{\xi_k \notin \bB_{r_n}(0) \} \right] \ton 1.
\end{align*}
The convergence to $1$ that we claim in the last step was established in~\eqref{eq:convergence of squares estimate of probability}.
Putting the results together we have
\begin{align*}
\P[P_n \text{ has exactly one critical point in } \bB_{r_n}(0)]
&=
\P\left[ \frac{P_n'}{P_n} \text{ has exactly one zero in } \bB_{r_n}(0) \right]\\
&\geq
\P[E_n^{(1)} \cap E_n^{(2)}] \ton 1.
\end{align*}
This completes the proof of~\eqref{eq:mainthm1}.

\medskip
\noindent
\textit{Proof of~\eqref{eq:mainthm2}.}
Recall that without restriction of generality we assume that $u_0=0$.
Also, from~\eqref{eq:z_n_w_def} we recall  the notation
\begin{equation}\label{eq:z_n_w_def2}
z_n(w)=  - \frac{1}{n f(0)} + \frac{w \sqrt{\log n}}{n^{\frac{3}{2}} (f(0))^2 }, \quad w\in\C.
\end{equation}
Observe that $z_n:\C\to\C$ is an affine function. Its inverse is denoted by $z_n^{-1}$. Recall that $\zeta_n$ is the zero of $P_n'$ in the disk $\bB_{r_n}(0)$ if it exists uniquely, and $\zeta_n$ is arbitrary otherwise. Our aim is to show that
\begin{equation}\label{eq:z_n_weak_need1}
z_n^{-1}(\zeta_n) \todistr N \sim \Normal_\C(0, \pi p(0)).
\end{equation}
Take some bounded, continuous function $\varphi: \C\to \R$. We need to show that
\begin{equation}\label{eq:z_n_weak_need2}
\lim_{n\to\infty} \E \varphi(z_n^{-1}(\zeta_n)) = \E \varphi(N).
\end{equation}
It suffices to assume that $\varphi\geq 0$ since in the general case we can write $\varphi= \varphi_+ - \varphi_-$ and use linearity.

Let $\rho>0$ and recall that $\bA_\rho$ denotes the Banach space of functions that are continuous on the closed disk $\bB_\rho(0)$ and analytic in its interior, endowed with the supremum norm.  Theorem~\ref{theo:funct.CLT} implies that
\begin{multline*}
\left(
\frac{1}{\sqrt{n \log n}} \left( \frac{P_n'(z_n(w))}{P_n(z_n(w))} -n f(0) - \frac{1}{z_n(w)} \right)
\right)_{|w| \le \rho}
\\
=
\left(
\frac{1}{\sqrt{n \log n}} \sum_{k=1}^n
\left( \frac{1}{z_n(w)-\xi_k} - f(0) \right)
\right)_{|w| \le \rho}
\toweak  (N)_{|w|\le \rho}
\end{multline*}
weakly on the space $\bA_\rho$. Here, $N \sim \Normal_\C(0, \pi p(0))$ is as above a complex normal random variable.
Using the Taylor expansion
\begin{align}
\frac{1}{z_n(w)} = -n f(0) - w \sqrt{n \log n} + o\left(\sqrt{n \log n}\right),
\quad \text{ as } n\to\infty,
\end{align}
that follows from~\eqref{eq:z_n_w_def2},  and applying Slutsky's lemma, we obtain
\begin{align}
\left(
\frac{1}{\sqrt{n \log n}} \frac{P_n'(z_n(w))}{P_n(z_n(w))} + w
\right)_{|w|\le \rho}
\toweak (N)_{|w|\le \rho} \label{eq:proofmainthm(b)_main_convergence}
\end{align}
weakly on $\bA_\rho$. Note that $f(\cdot)\mapsto f(\cdot)+\cdot$ defines a continuous mapping from $\bA_\rho$ to itself. Using the continuous mapping theorem we deduce that
\begin{align}
\left(
\frac{1}{\sqrt{n \log n}} \frac{P_n'(z_n(w))}{P_n(z_n(w))}
\right)_{|w|\le \rho}
\toweak (N - w)_{|w|\le \rho} \label{eq:proofmainthm(b)_main_convergence1}
\end{align}
weakly on $\bA_\rho$. The idea of what follows is quite simple. Assuming that $\rho$ is sufficiently large, the right-hand side has a unique zero at $N$. On the other hand, if there is a unique critical point $\zeta_n$ of $P_n$ in the disk $\bB_{r_n}(0)$ (which has probability converging to $1$ as $n\to\infty$), then the left-hand side has a zero at $z_n^{-1}(\zeta_n)$. Given the weak convergence~\eqref{eq:proofmainthm(b)_main_convergence1}, it is natural to conjecture the distributional convergence of the corresponding zeroes, that is $z_n^{-1}(\zeta_n)\to N$ in distribution. This would yield~\eqref{eq:z_n_weak_need1}.

In the following we shall justify the above  heuristics. First of all, recall our convention, see Remark~\ref{rem:convention}, that the function on the left-hand side of~\eqref{eq:proofmainthm(b)_main_convergence1} is considered on the event that its denominator has no zeroes. Strictly speaking, \eqref{eq:proofmainthm(b)_main_convergence1} means that
\begin{align}
\left(
\frac{1}{\sqrt{n \log n}} \frac{P_n'(z_n(w))}{P_n(z_n(w))} \ind_{G_n(\rho)}
\right)_{|w|\le \rho}
\toweak (N - w)_{|w|\le \rho} \label{eq:proofmainthm(b)_main_convergence2}
\end{align}
where we defined the random event
$$
G_n(\rho)   = \{\text{the disk $z_n(\bB_{\rho}(0))$ contains no zeroes of $P_n$}\}.
$$
Note that $\lim_{n\to\infty} \P[G_n(\rho)] = 1$, as we have already shown in the proof of~\eqref{eq:convergence of squares estimate of probability}.  Consider the functional $\Psi:\bA_\rho\to \R$ defined as follows:
$$
\Psi(f)
=
\begin{cases}
\varphi(z), &\text{ if $f$ has exactly one zero, denoted by $z$, in the interior of $\bB_{\rho/2}(0)$, }\\
0, &\text{ otherwise},
\end{cases}
$$
where $f\in \bA_\rho$. As always,  zeroes of analytic functions are counted with multiplicities. By the Hurwitz theorem on the zeroes of a perturbed analytic function~\cite[p.~152]{conway_book}, the functional $\Psi$ is continuous at every $f\in A$, where
\begin{equation}\label{eq:exceptional_set_A}
A=\{f\in \bA_\rho\colon f \text{ has no zeroes on the boundary of } \bB_{\rho/2}(0)\}.
\end{equation}
Note that the sample paths of the process $(N-w)_{|w|\le \rho}$ belong to $A$ with probability $1$. The continuous mapping theorem applied to~\eqref{eq:proofmainthm(b)_main_convergence2} yields
$$
\Psi\left(
\frac{1}{\sqrt{n \log n}} \frac{P_n'(z_n(\cdot))}{P_n(z_n(\cdot))} \ind_{G_n(\rho)}
\right)
\todistr
\Psi(N-\cdot).
$$
Since $|\Psi(f)|\leq \sup_{z\in\C} \varphi (z) <\infty$ for every $f\in\bA_\rho$, we can pass to expectations obtaining
\begin{equation}\label{eq:conv_expect_Psi}
\lim_{n\to\infty} \E \Psi\left(
\frac{1}{\sqrt{n \log n}} \frac{P_n'(z_n(\cdot))}{P_n(z_n(\cdot))} \ind_{G_n(\rho)}
\right)
=
\E \Psi(N-\cdot).
\end{equation}
By the definition of $\Psi$, we have
$$
\Psi(N-\cdot) = \varphi(N) \ind_{\{|N| < \rho/2\}}.
$$
In addition to the event $G_n(\rho)$ defined above consider the random events
\begin{align*}
F_n(\rho/2) &= \{\text{the disk $z_n(\bB_{\rho/2}(0))$ contains  exactly one zero of $P_n'$}\},\\
H_n &=\{\text{the disk $\bB_{r_n}(0)$ contains exactly one zero of  $P_n'$}\}.
\end{align*}
By~\eqref{eq:convergence of squares estimate of probability} and part (a) of Theorem~\ref{theo:criticalpoints}, the probability of $G_n(\rho)\cap H_n$ converges to $1$ as $n\to\infty$. In the following, we work on the event $G_n(\rho)\cap H_n$.  Observe that for sufficiently large $n$, we have $z_n(\bB_{\rho/2}(0))\subset \bB_{r_n}(0)$. Hence, on the event $F_n(\rho/2)\cap G_n(\rho)\cap H_n$, the only zero of $P_n'$ in $z_n(\bB_{\rho/2}(0))$ is $\zeta_n$. On the other hand, on the event  $(G_n(\rho)\cap H_n)\backslash F_n(\rho/2)$, there are no zeroes of $P_n'$ in $z_n(\bB_{\rho/2}(0))$. It follows that
$$
\Psi\left(
\frac{1}{\sqrt{n \log n}} \frac{P_n'(z_n(\cdot))}{P_n(z_n(\cdot))} \ind_{G_n(\rho)}
\right)
=
\varphi(z_n^{-1}(\zeta_n))
\ind_{F_n(\rho/2)}
\quad \text{ on the event } G_n(\rho)\cap H_n.
$$
Thus, we can write~\eqref{eq:conv_expect_Psi} in the following form:
$$
\lim_{n\to\infty}\E \left[ \varphi(z_n^{-1}(\zeta_n)) \ind_{F_n(\rho/2)} \right]
=
\E \left[\varphi(N) \ind_{\{|N| < \rho/2\}}\right].
$$
This looks almost like the desired statement~\eqref{eq:z_n_weak_need2} except that we still need to remove the indicator functions. This can be done as follows. Using that $\varphi$ is non-negative, we can write
$$
\liminf_{n\to\infty} \E \left[ \varphi(z_n^{-1}(\zeta_n))\right]
\geq
\lim_{n\to\infty} \E \left[ \varphi(z_n^{-1}(\zeta_n)) \ind_{F_n(\rho/2)}\right]
=
\E \left[\varphi(N) \ind_{\{|N| < \rho/2\}}\right].
$$
Since this holds for every $\rho>0$, we may let $\rho\to +\infty$, which yields the lower bound
$$
\liminf_{n\to\infty} \E \left[ \varphi(z_n^{-1}(\zeta_n))\right] \geq
\E \left[\varphi(N)\right].
$$
To prove the upper bound, observe that
$$
\E \left[ \varphi(z_n^{-1}(\zeta_n))\right]
\leq
\E \left[ \varphi(z_n^{-1}(\zeta_n))\ind_{F_n(\rho/2)}\right]
+
\|\varphi\|_{\infty}  \P[F_n^c(\rho/2)],
$$
where $\|\varphi\|_\infty = \sup_{z\in\C} |\varphi(z)|$. It follows that for every $\rho>0$,
$$
\limsup_{n\to\infty} \E \left[ \varphi(z_n^{-1}(\zeta_n))\right]
\leq
\E \left[\varphi(N) \ind_{\{|N| < \rho/2\}}\right]
+
\|\varphi\|_{\infty} \limsup_{n\to\infty} \P[F_n^c(\rho/2)],
$$
and to complete the proof of the upper bound it remains to show that
\begin{equation}\label{eq:limsup_F_n_compl}
\lim_{\rho\to +\infty} \limsup_{n\to\infty} \P[F_n^c(\rho/2)] = 0.
\end{equation}
By definition, $F_n^c(\rho/2)$ is the union of the following two events:
\begin{align*}
F_n^{(0)}(\rho/2) &= \{\text{the disk $z_n(\bB_{\rho/2}(0))$ contains no zeroes of $P_n'$}\},\\
F_n^{(\geq 2)}(\rho/2) &= \{\text{the disk $z_n(\bB_{\rho/2}(0))$ contains at least two zeroes of $P_n'$}\}.
\end{align*}
Since for sufficiently large $n$, we have  $z_n(\bB_{\rho/2}(0))\subset \bB_{r_n}(0)$ and hence the event $F_n^{(\geq 2)}(\rho/2)$ is contained in the complement of $H_n$, we have
\begin{equation}\label{eq:F_n_geq_2}
\lim_{n\to\infty} \P[F_n^{(\geq 2)}(\rho/2)] = 0.
\end{equation}
To estimate the probability of $F_n^{(0)}(\rho/2)$, consider the following functional $\Psi_0:\bA_\rho \to \R$:
$$
\Psi_0(f)= \ind_{\{f \text{ has no zeroes in } \bB_{\rho/2}(0)\}}.
$$
Clearly, $\Psi_0$ is continuous on $A$ defined in~\eqref{eq:exceptional_set_A}, and hence the same argumentation as in~\eqref{eq:conv_expect_Psi} leads to
$$
\lim_{n\to\infty} \E \Psi_0\left(
\frac{1}{\sqrt{n \log n}} \frac{P_n'(z_n(\cdot))}{P_n(z_n(\cdot))} \ind_{G_n(\rho)}
\right)
=
\E \Psi_0(N-\cdot) = \P[|N| > \rho/2].
$$
By definition of $\Psi_0$, this yields
$$
\limsup_{n\to\infty} \P[F_n^{(0)}(\rho/2)] \leq \lim_{n\to\infty} \E \Psi_0\left(
\frac{1}{\sqrt{n \log n}} \frac{P_n'(z_n(\cdot))}{P_n(z_n(\cdot))} \ind_{G_n(\rho)}
\right)
=
\P[|N| > \rho/2].
$$
Letting $\rho\to +\infty$ and recalling~\eqref{eq:F_n_geq_2}, we arrive at~\eqref{eq:limsup_F_n_compl}, thus proving~\eqref{eq:z_n_weak_need2}.
The proof of~\eqref{eq:mainthm2} and of Theorem~\ref{theo:criticalpoints} is complete.
\hfill $\Box$

\subsection{Proof of Theorem \ref{theo:criticalpoints_random_u0}}
Without restriction of generality, let $p>0$ on $\mathscr D$ (otherwise, we can replace $\mathscr D$ by $\mathscr D\backslash \{p=0\}$).
The conditional distribution of $Q_n$ given that $\xi_0 = u_0$ can be identified with that of $P_n$. Part (a) of Theorem~\ref{theo:criticalpoints} implies
\begin{align}
\lim_{n \to \infty} \P \left[\text{there is unique }  \zeta \in \bB_{r_n}(\xi_0) \text{ such that }  Q_n'(\zeta)=0 \middle| \xi_0 = u_0 \right] = 1
\end{align}
for all $u_0\in\mathscr D$ such that  $f(u_0)\neq 0$. Recall that this condition is violated only on a set of Lebesgue measure $0$ of $u_0$'s.
Hence, by the total probability formula and the  dominated convergence theorem
\begin{align*}
& \P [\text{there is unique }  \zeta \in \bB_{r_n}(\xi_0) \text{ such that }  Q_n'(\zeta)=0] \\
&=\int_{\mathscr D} \P \left[\text{there is unique }  \zeta \in \bB_{r_n}(\xi_0) \text{ such that }  Q_n'(\zeta)=0 \middle| \xi_0 = u_0 \right]  p(u_0) \dd u_0\\
&\ton  \int_{\mathscr D} p(u_0) \dd u_0 =1.
\end{align*}
To prove part (b) of the theorem, recall that a sequence of complex-valued random variables $(X_n)_{n\in\N}$ converges to a random variable $X$ in distribution iff $\lim_{n\to\infty} \E \varphi(X_n) = \E \varphi(X)$ for every bounded continuous function $\varphi: \C\to \R$. So, let $\varphi$ be such function. Then, part~(b) of Theorem~\ref{theo:criticalpoints} and the definition of weak convergence yield
\begin{multline*}
\lim_{n\to\infty}
\E_{Q_n} \left[ \varphi\left(\frac{f^2 (\xi_0)}{\sqrt{\pi p(\xi_0)}} \sqrt{\frac{n}{\log n}}   \left( n \left( \zeta_n-\xi_0 \right) +\frac{1}{f \left( \xi_0 \right)}\right)\right) \middle | \xi_0 = u_0\right]
\\
=
\lim_{n\to\infty}
\E_{P_n} \left[ \varphi\left(\frac{f^2 (u_0)}{\sqrt{\pi p(u_0)}} \sqrt{\frac{n}{\log n}}   \left( n \left( \zeta_n-u_0 \right) +\frac{1}{f \left(u_0 \right)}\right)\right)\right]
=
\E \varphi(N),
\end{multline*}
for all $u_0\in\mathscr D$ such that  $f(u_0)\neq 0$. Here, $N\sim \mathcal N_\C(0,1)$ is a standard complex normal variable.
Using this together with the dominated convergence theorem, we obtain
\begin{multline*}
\E_{Q_n} \left[ \varphi\left(\frac{f^2 (\xi_0)}{\sqrt{\pi p(\xi_0)}} \sqrt{\frac{n}{\log n}}   \left( n \left( \zeta_n-\xi_0 \right) +\frac{1}{f \left( \xi_0 \right)}\right)\right)\right]
\\=
\int_{\mathscr D} \E_{Q_n} \left[ \varphi\left(\frac{f^2 (\xi_0)}{\sqrt{\pi p(\xi_0)}} \sqrt{\frac{n}{\log n}}   \left( n \left( \zeta_n-\xi_0 \right) +\frac{1}{f \left( \xi_0 \right)}\right)\right) \middle | \xi_0 = u_0 \right]
p(u_0) \dd u_0
\\
\ton
\int_{\mathscr D} \E \varphi(N)p(u_0) \dd u_0
=
\E \varphi(N).
\end{multline*}
Since this holds for every bounded continuous function $\varphi$, the proof is complete.
\hfill $\Box$


\bibliography{critical_points_bib}
\bibliographystyle{plainnat}

\end{document}